\documentclass[12pt]{amsart}
\usepackage{verbatim,amssymb,latexsym,amscd}
\usepackage[all]{xy}

\usepackage[colorlinks,linkcolor=blue,citecolor=blue,urlcolor=red]{hyperref}
\newcommand{\ie}{{\it i.e. }}
\newcommand{\cf}{{\it cf. }}
\newcommand{\eg}{{\it e.g. }}
\newcommand{\loccit}{{\it loc. cit. }}
\newcommand{\resp}{{\it resp. }}
\newcommand{\un}{\mathbf{1}}
\newcommand{\A}{\mathbf{A}}

\newcommand{\C}{\mathbf{C}}

\renewcommand{\P}{\mathbf{P}}
\newcommand{\Q}{\mathbf{Q}}

\newcommand{\Z}{\mathbf{Z}}
\newcommand{\sA}{\mathcal{A}}
\newcommand{\sB}{\mathcal{B}}
\newcommand{\sC}{\mathcal{C}}
\newcommand{\sD}{\mathcal{D}}

\newcommand{\sF}{\mathcal{F}}
\newcommand{\sG}{\mathcal{G}}

\newcommand{\sI}{\mathcal{I}}
\newcommand{\sJ}{\mathcal{J}}
\newcommand{\sK}{\mathcal{K}}
\newcommand{\sL}{\mathcal{L}}
\newcommand{\sM}{\mathcal{M}}

\newcommand{\sR}{\mathcal{R}}

\newcommand{\sZ}{\mathcal{Z}}

\newcommand{\bL}{\mathbb{L}}
\newcommand{\bT}{\mathbb{T}}

\newcommand{\Spec}{\operatorname{Spec}}
\newcommand{\Ker}{\operatorname{Ker}}
\newcommand{\Coker}{\operatorname{Coker}}
\newcommand{\IM}{\operatorname{Im}}
\newcommand{\place}{\operatorname{\bf place}}
\newcommand{\field}{\operatorname{\bf field}}

\newcommand{\Coprod}{{\coprod}}

\newcommand{\Mor}{\operatorname{Mor}}
\newcommand{\Hom}{\operatorname{Hom}}
\newcommand{\End}{\operatorname{End}}
\newcommand{\Mod}{\text{\rm Mod--}}

\newcommand{\Albv}{\operatorname{Alb}}
\newcommand{\NS}{\operatorname{NS}}
\newcommand{\Pic}{\operatorname{Pic}}
\newcommand{\Griff}{\operatorname{Griff}}

\newcommand{\cl}{{\operatorname{cl}}}

\newcommand{\Ext}{\operatorname{Ext}}
\newcommand{\AbS}{\operatorname{\bf AbS}}
\newcommand{\Ab}{\operatorname{\bf Ab}}
\newcommand{\SAbS}{\operatorname{\bf SAbS}}
\newcommand{\SAb}{\operatorname{\bf SAb}}
\newcommand{\Lat}{\operatorname{\bf Lat}}
\newcommand{\Mot}{\operatorname{\bf Mot}}

\newcommand{\car}{\operatorname{char}}
\newcommand{\tors}{{\operatorname{tors}}}
\newcommand{\rat}{{\operatorname{rat}}}
\newcommand{\alg}{{\operatorname{alg}}}
\renewcommand{\hom}{{\operatorname{hom}}}
\newcommand{\num}{{\operatorname{num}}}
\newcommand{\AJ}{{\operatorname{AJ}}}

\newcommand{\tnil}{{\operatorname{tnil}}}
\newcommand{\CM}{\operatorname{\bf CM}}

\newcommand{\Sm}{\operatorname{\bf Sm}}

\newcommand{\proper}{{\operatorname{prop}}}
\newcommand{\qp}{{\operatorname{qp}}}

\newcommand{\proj}{{\operatorname{proj}}}

\newcommand{\Cor}{\operatorname{\bf Cor}}

\newcommand{\Chow}{\operatorname{\bf Chow}}

\newcommand{\rep}{{\text{\rm rep}}}
\newcommand{\et}{{\text{\rm \'et}}}
\renewcommand{\o}{{\text{\rm o}}}
\newcommand{\op}{{\text{\rm op}}}
\renewcommand{\b}{{\text{\rm b}}}
\newcommand{\eff}{{\text{\rm eff}}}

\newcommand{\ttto}{\dashrightarrow}
\newcommand{\tto}{\rightsquigarrow}

\newcommand{\rr}{\operatornamewithlimits{\rightrightarrows}\limits}
\newcommand{\by}[1]{\overset{#1}{\longrightarrow}}

\newcommand{\iso}{\by{\sim}}

\newcommand{\inj}{\hookrightarrow}

\newcommand{\surj}{\rightarrow\!\!\!\!\!\rightarrow}
 
\newcommand{\colim}{\varinjlim}
\renewcommand{\lim}{\varprojlim}

\renewcommand{\qed}{\hfill $\Box$\medskip}

\renewcommand{\phi}{\varphi}
\renewcommand{\epsilon}{\varepsilon}

\newcounter{spec}
\newenvironment{thlist}{\begin{list}{\rm{(\roman{spec})}}%
{\usecounter{spec}\labelwidth=20pt\itemindent=0pt\labelsep=10pt}}%
{\end{list}}%

\swapnumbers

\newtheorem{thm}{Theorem}[subsection]
\newtheorem{lemma}[thm]{Lemma}
\newtheorem{prop}[thm]{Proposition}
\newtheorem{cor}[thm]{Corollary}
\newtheorem{conj}[thm]{Conjecture}

\theoremstyle{definition}

\newtheorem{defn}[thm]{Definition}

\newtheorem{rk}[thm]{Remark}
\newtheorem{rks}[thm]{Remarks}

\newtheorem{ex}[thm]{Example}
\newtheorem{exs}[thm]{Examples}

\numberwithin{equation}{section}

\entrymodifiers={!!<0pt,0.7ex>+}

\begin{document}
\title{Birational motives, I: pure birational motives}
\author{Bruno Kahn}
\address{IMJ-PRG\\ Case 247\\4 place Jussieu\\
75252 Paris Cedex 05\\France}
\email{bruno.kahn@imj-prg.fr}
\author{R. Sujatha}
\address{University of British Columbia\\Vancouver, BC
V6T1Z2\\Canada}
\email{sujatha@math.ubc.ca}
\date{October 28, 2015}
\thanks{The first author acknowledges the support of Agence Nationale de la Recherche (ANR) under reference ANR-12-BL01-0005 and the second author  that of NSERC Grant 402071/2011. Both authors acknowledge the support of CEFIPRA project 2501-1.}
\begin{abstract}
We define a category of pure birational motives over a field, depending on the choice of an adequate equivalence relation on algebraic cycles. It is obtained by ``killing'' the Lefschetz motive in the corresponding category of effective motives. For rational equivalence, it encompasses Bloch's decomposition of the diagonal. We study the induced Chow-K\"unneth decompositions in this category, and establish relationships with Rost's cycle modules and the Albanese functor for smooth projective varieties.
\end{abstract}
\maketitle

\tableofcontents

\section*{Introduction}

In the preprint \cite{birat}, we toyed with birational ideas in three 
areas of algebraic geometry: plain varieties, pure motives in the sense 
of Grothendieck, and triangulated motives in the sense of Voevodsky. 
These three themes are finally treated separately in revised versions. 
The first one is the object of \cite{Birat}; the second one is the 
object of the present paper; the third one is the object of \cite{birat-tri}.

We work over a field $F$. Recall that we introduced in \cite{Birat} two
``birational" categories. The first, $\place(F)$, has for objects the function
fields over $F$ and for morphisms the $F$-places. The second one is the
Gabriel-Zisman localisation of the category $\Sm(F)$ of smooth $F$-varieties
obtained by inverting birational morphisms \cite[Ch. 1]{gz}: we denoted this category by
$S_b^{-1}\Sm(F)$.

We may also invert stable birational morphisms: those which are dominant and
induce a purely transcendental extension of function fields, and invert the
corresponding morphisms in $\place(F)$. We denote the sets of such morphisms by
$S_r$.

In order to simplify the exposition, let us assume that $F$ is of 
characteristic
$0$. Then the main results of \cite{Birat} and its predecessor
\cite{localisation} can be summarised in a diagram 
\[\begin{CD}
\place(F)^\op@>>> S_b^{-1} \Sm^\proj(F)@>\sim>> S_b^{-1} \Sm(F)\\
@V{}VV @V{\wr}VV @V{\wr}VV \\
S_r^{-1}\place(F)^\op@>>> S_r^{-1} \Sm^\proj(F)@>\sim>> S_r^{-1} \Sm(F)
\end{CD}\]
where $\Sm^\proj(F)$ is the full subcategory of smooth projective 
varieties and
the symbols $\sim$ denote equivalences of categories: see
\cite[Prop. 8.5]{localisation} and \cite[Th. 1.7.2 and 
4.2.4]{Birat}. 

Moreover, if
$X$ is smooth and $Y$ is smooth proper, then $\Hom(X,Y)=Y(F(X))/R$ in
$S_b^{-1}\Sm(F)$, where $R$ is R-equivalence \cite[Th. 6.6.3]{Birat}.

In this paper, we consider the effect of inverting birational morphisms in
categories of \emph{effective pure motives}. For simplicity, let us 
still assume
$\car F=0$, and consider only the category of effective Chow motives
$\Chow^\eff(F)$, defined by using algebraic cycles modulo rational 
equivalence. 
The graph functor then induces a commutative square (compare \eqref{eq5.1})
\[\begin{CD}
S_b^{-1} \Sm^\proj(F)@>>> S_b^{-1} \Chow^\eff(F)\\
@V{\wr}VV @V{}VV \\
S_r^{-1} \Sm^\proj(F)@>>> S_r^{-1} \Chow^\eff(F).
\end{CD}\]

One can expect that the right vertical functor is an equivalence of 
categories,
and indeed this is not difficult to prove (Corollary \ref{c7.1} b)). 
But we have
two other descriptions of this category of ``birational motives":

\begin{itemize}
\item The functor $\Chow^\eff(F)\to S_b^{-1}\Chow^\eff(F)$ is full, 
and its kernel is the ideal $\sL_\rat$ of morphisms which factor 
through some object of
the form $M\otimes\bL$, where $\bL$ is the \emph{Lefschetz motive} (ibid.).
\item If $X,Y$ are smooth projective varieties, then 
$\sL_\rat(h(X),h(Y))$
coincides with the group of Chow correspondences represented by 
algebraic cycles on $X\times Y$ whose irreducible components are not 
dominant over $X$ (Theorem \ref{l5.1}).
\end{itemize}

As a consequence, the group of morphisms from $h(X)$ to $h(Y)$
in $S_b^{-1}
\Chow^\eff(F)$ is isomorphic to $CH_0(Y_{F(X)})$. Given the similar
description of Hom sets in $S_b^{-1}\Sm^\proj(F)$ recalled above, this 
places the classical map
\[Y(F(X))/R\to CH_0(Y_{F(X)})\]
in a categorical context. 

Note that, by \cite[Th. 8.5.1 b)]{Birat}, if $X\simeq \Spec F$ in $S_b^{-1} \Sm$ then $X$ must be rationally connected; on the other hand, there are surfaces of general type with trivial birational motive, see Remarks \ref{r3.1} 1) and 3). So the birational motive of a smooth projective variety detects much less geometry than its class in $S_b^{-1} \Sm$, but on the other hand it is much more computable.

This paper is organised as follows. In Section \ref{s1} we review pure motives.
In Section \ref{s2} we study pure birational motives, in greater generality
than outlined in this introduction. In particular, many results are valid for
other adequate equivalence relations than rational equivalence, see 
\S\ref{s2.3}; moreover, most
results extend to characteristic $p$ if $p$ is invertible in the ring of coefficients, 
by using the de Jong-Gabber alteration theorem \cite{gabber}, see Theorem 
\ref{l5.1}.

Section \ref{sexamples} consists of examples. We study varieties whose birational motive is trivial, in the line of the remarks above. We 
also study the Chow-K\"unneth decomposition in the category of 
birational motives, special attention being devoted to the case of 
complete intersections.

Let $\Chow^\o(F)$ denote the pseudo-abelian envelope of $S_b^{-1} \Chow^\eff(F)$. In Section
\ref{s3}, we examine two questions: the existence of a  right
adjoint to the projection functor $\Chow^\eff(F)\to \Chow^\o(F)$ (and
similarly for more general adequate equivalences), and whether pseudo-abelian completion is
really necessary. It turns out that the
answer to the first question is negative (Theorems \ref{t3.1} and \ref{t3.2}; this is related to the nontriviality of the Griffiths group for some $3$-folds)
and the answer to the second question is positive with rational coefficients  under a
nilpotence conjecture (Conjecture \ref{c3.1}). We can get  an unconditional positive
answer to the second question if we restrict to a suitable type of motives (Proposition
\ref{p3.1} and Example \ref{e3.1}).

In Section \ref{s4}, we define a functor $S_r^{-1}\field(F)^\op\to
S_r^{-1}\Chow^\eff(F,\Q)$ in characteristic $p$, using de Jong's 
theorem again. Here $\field(F)$ denotes the subcategory of $\place(F)$ 
with the same objects but morphisms restricted to field extensions 
(Proposition \ref{ph0}). 

We end this paper by relating the previous constructions to more
classical objects. In Section \ref{s.6} we relate birational motives to cycle cohomology \cite{rost}, expanding a bit on previous results by Rost and  Merkurjev \cite{merk,merk2}. In Section \ref{s5}, we define a tensor additive category
$\AbS(F)$ of \emph{locally abelian schemes}, whose objects are those $F$-group
schemes that are extensions of a lattice (\ie locally isomorphic for the
\'etale topology to a free finitely generated abelian group) by an
abelian variety. We then show in Section \ref{s6} that the classical construction
of the Albanese variety of a smooth projective variety extends to a tensor
functor
\[\Albv:\Chow^\o(F)\to \AbS(F)\]
which becomes full and essentially surjective after tensoring morphisms
with $\Q$ (Proposition \ref{p6.1}). So, one could say that $\AbS(F)$ is the
\emph{representable part} of $\Chow^\o(F)$. We also show that, after
tensoring with $\Q$, $\Albv$ has a right adjoint which
identifies $\AbS(F)\otimes\Q$ with the thick subcategory of
$\Chow^\o(F)\otimes\Q$ generated by motives of varieties of dimension
$\le 1$.
\medskip

Some results of the preliminary version \cite{birat} of this work were used in
other papers, namely \cite{kmp} and \cite{FqX}, and we occasionally refer to
these papers to ease the exposition. Here is a correspondence guide between the results from \cite{birat} used in these papers and those in the present version:

\begin{itemize}
\item In \cite{FqX}, Lemma 7.2 uses \cite[Lemmas 5.3 and 5.4]{birat}, which
correspond to Proposition \ref{l4.2'} and Theorem \ref{l5.1} of the present
paper. The reader will verify that the proofs of Proposition \ref{l4.2'} and
Theorem \ref{l5.1} are the same as those of \cite[Lemmas 5.3 and 5.4]{birat},
mutatis mutandis, and do not use any result from \cite{FqX}.
\item In \cite{kmp}, Lemma 7.5.3 uses the same references: the same comment as
above applies. Moreover, \cite[9.5]{birat} is used on pp. 174--175 of \cite{kmp}:
this result is now Theorem \ref{p6.2}. Again, its proof is identical to the
one in the preliminary version and does not use results from \cite{kmp}. 
\end{itemize}

The idea of considering birational Chow correspondences, that yield 
here a 
category in which $\Hom([X],[Y])=CH_0(Y_{F(X)})$ for two smooth 
projective varieties $X,Y$, goes back to S. Bloch's method of 
``decomposition of the diagonal" in \cite[App. to Lecture 1]{bloch} (see also Bloch-Srinivas \cite{bs}). 
He attributes the idea of considering the 
generic point of a smooth projective variety $X$ as a $0$-cycle over its function field to 
Colliot-Th\'el\`ene: here, this corresponds to the identity endomorphism of $h^\o(X)\in \Chow^\o(F)$. We 
realised the connection with Bloch's ideas after reading H. Esnault's 
article 
\cite{esnault}, and this led to another proof of her theorem by the 
present birational techniques in \cite{FqX}. M. Rost has considered 
this category independently \cite{merk}: this was pointed out to us by 
N. Karpenko.

\subsection*{Acknowledgements} We thank Andr\'e, Barbieri-Viale, Bondarko, Col\-liot-Th\'el\`ene, Fakhruddin, Karpenko, Maltsiniotis, Srinivas and Voisin for helpful exchanges. Special thanks to the referee, whose enthusiastic report prompted us to considerably expand Subsection \ref{s3.1} of this paper.

\section{Review of pure motives}\label{s1}

In this section, we recall the definition of categories of pure motives in a way
which is suited to our needs. A slight variance to the usual exposition is the
notion of \emph{adequate pair} which is a little more precise than the notion of
adequate equivalence relation (it explicitly takes the coefficients into
account).

We adopt the covariant convention, for future comparison with Voevodsky's
triangulated categories of motives: here, the functor which sends a smooth
projective variety to its motive is covariant. For a dictionary between the
covariant and contravariant conventions, the reader may refer to
\cite[7.1.2]{kmp}.

\subsection{Adequate pairs}

We give ourselves:

\begin{itemize}
\item a commutative ring of coefficients $A$;
\item an adequate equivalence relation $\sim$ on algebraic cycles
with coefficients in $A$ \cite{samuel}.
\end{itemize}

We refer to $(A,\sim)$ as an \emph{adequate pair}. Classical examples for $\sim$ are $\rat$
(rational equivalence),
$\alg$ (algebraic equivalence), $\num$ (numerical equivalence), $\sim_H$
(homological equivalence relative to a fixed Weil cohomology theory
$H$). A less classical example is Voevodsky's smash-nilpotence $\tnil$
\cite{voenil}, see \cite[Ex. 7.4.3]{ak} (a cycle $\alpha$ is
smash-nilpotent if $\alpha^{\otimes n}\sim_\rat 0$ for some $n>0$).
We then have a notion of domination
$(A,\sim)\ge (A,\sim')$ if $\sim$ is finer than $\sim'$ (\ie the groups
of cycles modulo $\sim$ surjects onto the one for $\sim'$). It is well-known
that $(A,\rat)\ge (A,\sim)$ for any $\sim$ (\cf \cite[Ex.
1.7.5]{fulton}), and that $(A,\sim)\ge (A,\num_A)$ if $A$ is a
field.

Since the issue of coefficients is sometimes confusing, the following
remarks may be helpful. Given a pair $(A,\sim)$ and a commutative
$A$-algebra $B$, we get a new pair $B\otimes_A (A,\sim)$
by tensoring algebraic cycles with $B$: for example,
$(A,\sim)=A\otimes_\Z (\Z,\sim)$ for $\sim=\rat,\alg$ or $\tnil$ by
definition. On the other hand, given a pair $(B,\sim)$ and a ring
homomorphism $A\to B$ we get a ``restriction of scalars" pair
$(A,\sim_{|A})$ by considering cycles with coefficients in $A$ which
become $\sim 0$ after tensoring with $B$: for example, if $H$ is a
Weil cohomology theory with coefficients in $K$, this applies to
any ring homomorphism $A\to K$. Obviously $B\otimes_A
(A,\sim_{|A})\ge (B,\sim)$, but this need not be an equality in general.

In the case of numerical equivalence (a cycle with coefficients in $A$ is
numerically equivalent to $0$ if the degree of its intersection with any
cycle of complementary dimension in good position is $0$), we have
$B\otimes_A (A,\num_A)\ge (B,\num_B)$, with equality if $B$ is flat
over $A$.

Given a pair $(A,\sim)$, to any smooth projective $F$-variety $X$ we may
associate for each integer $n\ge 0$ its group of cycles of
codimension $n$ with coefficients in $A$ modulo
$\sim$, that will be denoted by $\sZ_\sim^n(X,A)$. If $X$ has pure
dimension $d$, we also write this group $\sZ_{d-n}^\sim(X,A)$.

\subsection{Smooth projective varieties, connected and nonconnected} In
\cite{Birat} we were only considering (connected) varieties over $F$.
Classically, pure motives are defined using not necessarily connected smooth
projective varieties. One could base the treatment on connected smooth
varieties, but this would introduce problems with the tensor product, since a
product of connected varieties need not be connected in general (\eg if neither
of them is geometrically connected). Thus we prefer to use here:

\begin{defn}\label{d1.3} We write $\Sm_\Coprod(F)$ for the category of smooth
separated schemes of finite type over $F$. For $\% \in \{\proper,\qp,\proj\}$,
we write
$\Sm^\%_\Coprod(F)$ for the full subcategory of $\Sm_\Coprod(F)$ consisting of
proper, quasi-projective or projective varieties.
\end{defn} 

Unlike their counterparts considered in \cite{Birat}, these categories enjoy
finite products and coproducts. 

The following lemma is clear.

\begin{lemma}\label{l1.3} The categories considered in Definition \ref{d1.3} are
the ``finite coproduct envelopes" of those considered in \cite{Birat}, in
the sense of \cite[Prop. 6.1]{localisation}.
\end{lemma}

\subsection{Review of correspondences}\label{s1.3}

We associate to two smooth projective varieties $X,Y$ the
group $\sZ^{\dim Y}_\sim(X\times Y,A)$ of correspondences from $X$ to $Y$
relative to $(A,\sim)$. The composition of correspondences is defined
as follows\footnote{We follow here the convention of
Voevodsky in \protect\cite{voetri}. It is also the one used by Fulton
\protect\cite[\S 16]{fulton}. See \protect\cite[7.1.2]{kmp}.}: if
$X,Y,Z$ are smooth projective and $(\alpha,\beta)\in \sZ^{\dim
Y}_\sim(X\times Y,A)\times \sZ^{\dim Z}_\sim(Y\times Z,A)$, then
\[\beta\circ\alpha = (p_{XZ})_*(p_{XY}^*\alpha\cdot p_{YZ}^*\beta)\]
where $p_{XY}, p_{YZ}$ and $p_{XZ}$ denote the partial projections from
$X\times Y\times Z$ onto two-fold factors.

We then get an $A$-linear tensor (\ie symmetric monoidal) category
$\Cor_\sim(F,A)$. The graph map defines a
\emph{covariant} functor
\begin{align}
\Sm_\Coprod^\proj(F)&\to \Cor_\sim(F,A)\label{eq7.1}\\
X&\mapsto [X]\notag
\end{align}
so that $[X\coprod Y] = [X]\oplus [Y]$, and $[X\times Y] = [X]\otimes [Y]$ for
the tensor structure. The unit object is $\un=[\Spec F]$.

If $f:X\to Y$ is a morphism of smooth varieties, let $\Gamma_f$ denote its graph and
$[\Gamma_f]$ denote the class of $\Gamma_f$ in $\sZ^{\dim
Y}_\sim(X\times Y)$. We write $f_*$ for the correspondence
$[\Gamma_f]:[X]\to [Y]$ (the image of $f$ under the functor
\eqref{eq7.1}). Note that if
$f:X\to Y$ and $g:Y\to Z$ are two morphisms of smooth projective varieties, then the cycles
$\Gamma_f\times Z$ and $X\times \Gamma_g$ on $X\times Y\times Z$
intersect properly, so that $g_*\circ f_*$ is well-defined as a
cycle and not just as an equivalence class of cycles; the equation
$g_*\circ f_*=(g\circ f)_*$ is an equality of cycles. (This is a very special case of the composition of finite correspondences, cf. \cite[Lemma 1.7]{mvw}.)

\subsection{The correspondence attached to a rational map}\label{s1.4} We first define rational maps between not necessarily
connected smooth varieties $X,Y$ in the obvious way: it is a morphism from a
suitable \emph{dense} open subset of $X$ to $Y$. Like morphisms, rational maps
split as disjoint unions of ``connected" rational maps. A rational map $f$ is
\emph{dominant} if all its connected components are dominant and if the image of $f$ meets all connected components of $Y$.

Let $f:X\ttto Y$ be a rational
map between two smooth projective varieties $X,Y$. To $f$ we associated in \cite[\S 6.3]{Birat} a morphism in the category $S_b^{-1}\Sm$. In the case of Chow motives, we can do better: define
the correspondence
$f_*:[X]\to [Y]$ in $\Cor_\sim(F,A)$,  as the closure of the graph
of $f$ inside $X\times Y$. The formula $g_*\circ f_*=(g\circ f)_*$ need not be
valid in general, even if $g\circ f$ is defined (but see Proposition \ref{l7.3}
below). Yet we have:

\begin{lemma}\label{l1.1} Let $X\overset{f}{\ttto} Y\by{g} Z$ be a diagram of
smooth projective varieties, where $f$ is a rational map and $g$ is a morphism.
Then we have an equality of cycles
\[g_*\circ f_*=(g\circ f)_*\]
in $\sZ^{\dim Z}(X\times Z)$.
\end{lemma}

\begin{proof} Let $U$ be an open subset of $X$ on which $f$, hence $g\circ f$, is defined. As explained in \ref{s1.3}, we have an equality of reduced closed subschemes $\Gamma_{g\circ f}= p_{UZ}(\Gamma_f\times Z\cap X\times \Gamma_g)$.  Since $Y$ is proper,  $p_{UZ}(\Gamma_f\times Z\cap X\times \Gamma_g)$ is dense in $p_{XZ}(\bar\Gamma_f\times Z\cap X\times \Gamma_g)=g_*\circ f_*$, hence the conclusion.
\end{proof}

\enlargethispage*{30pt}

\subsection{Effective pure motives}

We now define as usual the category of effective pure motives
$\Mot^\eff_\sim(F,A)$ relative to $(A,\sim)$ as the pseudo-abelian
envelope of $\Cor_\sim(F,A)$. We denote the composition of
\eqref{eq7.1} with the pseudo-abelianisation functor by $h_\sim$. If $\sim=\rat$, we
usually abbreviate $h_\sim$ to $h$. 

In $\Mot^\eff_\sim(F,A)$ we have
\begin{itemize}
\item $h_\sim(\Spec F)=\un$ (the unit object for the tensor
structure)
\item $h_\sim(\P^1)=\un\oplus \bL $ where $\bL $ is the \emph{Lefschetz
motive}.
\end{itemize}

If $n\ge 0$, we write $M(n)$ for the motive
$M\otimes \bL^{\otimes n}$ (beware that the ``standard" notation is
$M(-n)$!)

We then have the formula, for two smooth projective $X,Y$ and integers
$p,q\ge 0$
\begin{equation}\label{eq7.3}
\Mot^\eff_\sim(F,A)(h_\sim(X)(p),h_\sim(Y)(q))=\sZ_\sim^{\dim
Y+q-p}(X\times Y).
\end{equation}

In particular, the endofunctor $-\otimes \bL $ of $\Mot^\eff_\sim(F,A)$ is
fully faithful. 

If $f:X\to Y$ is a morphism, then the
correspondence  $[{}^t\Gamma_f]\in \sZ^{\dim Y}(Y\times X)$ obtained by
the ``switch" defines a morphism $f^*:h_\sim(Y)(\dim X)\to h_\sim(X)(\dim
Y)$, \ie from $h_\sim(Y)$ to $h_\sim(X)(\dim Y-\dim X)$ or from
$h_\sim(Y)(\dim X-\dim Y)$ to $h_\sim(X)$ according to the sign of $\dim
X-\dim Y$. In particular, if $f$ has relative dimension $0$ then $f^*$ maps $h_\sim(Y)$ to
$h_\sim(X)$.  We similarly define $f^*$ for a rational map $f$.

We recall the well-known

\begin{lemma}\label{l7.1} Suppose that $f$ is generically finite of
degree
$d$. Then $f_*\circ f^*=d 1_Y$.
\end{lemma}

\begin{proof} It suffices to prove this for the action on cycles, and then
the lemma follows by Manin's identity principle \cite[\S 2]{scholl}. Let $\alpha\in
\sZ^*_\sim(Y,A)$. By the projection formula,
\[f_*f^*(\alpha)=\alpha\cdot f_*(1).\]

But $f_*(1)\in \sZ^0_\sim(Y,A)$ may be computed after restriction
to any open subset $U$ of $X$ and for $U$ small enough it is clear
that $f_*(1)=d$.\end{proof}

\subsection{Pure motives} The category $\Mot_\sim(F,A)$ is now obtained from
$\Mot_\sim^\eff(F,A)$ by inverting the endofunctor $-\otimes\bL$, \ie adjoining a
$\otimes$-quasi-inverse $\bT$ of $\bL$ (the Tate motive) to $\Mot_\sim^\eff(F,A)$. The
resulting category is rigid and the functor $\Mot_\sim^\eff(F,A)\to \Mot_\sim(F,A)$ is fully
faithful; we refer to \cite{scholl} for details. We still write $h_\sim(X)$ for the image of
$h_\sim(X)$ in $\Mot_\sim(F,A)$.

\subsection{Pure motives and purely inseparable extensions} This subsection will be needed for the proof of Theorem \ref{rhomsim} below. It shows that extending scalars along a purely inseparable extension is harmless as long as the exponential characteristic is inverted.

\begin{lemma}\label{l4.1a} Let $f:X\to Y$ be a finite, flat and radicial  morphism \cite[Def. 3.7.2]{EGAI} between
smooth projective  $F$-varieties. Let $(A,\sim)$ be an adequate pair, with $p$ invertible
in $A$ (where $p$ is the exponential characteristic of $F$). Then\\ 
a) $f_*:\sZ_*^\sim(X,A)\to \sZ_*^\sim(Y,A)$ is an isomorphism.\\ 
b) $f_*:h(X)\to h(Y)$ is an isomorphism in $\Cor_\sim(F,A)$.
\end{lemma}

\begin{proof} Let $p^n$ be the generic  degree of $f$. We have $f_* f^*=f^*f_* = p^n$ (on the level of  algebraic cycles), hence a). b) follows by
Manin's identity principle (Yoneda lemma).
\end{proof}

\begin{prop}\label{l3.1a} Let $K/F$ be a purely inseparable extension. Then, for any 
adequate pair $(A,\sim)$ as in Lemma \ref{l4.1a}, the extension of scalars functors
\begin{align*}
\Cor_\sim(F,A)&\to \Cor_\sim(K,A)\\
\Mot_\sim^\eff(F,A)&\to \Mot_\sim^\eff(K,A)\\
\Mot_\sim(F,A)&\to \Mot_\sim(K,A)
\end{align*}
are equivalences of categories.
\end{prop}

\begin{proof} It suffices to show this for the first functor. Let
$X,Y$ be two smooth projective $F$-varieties. Then, for any finite sub-extension $L/F$ of
$K/F$, the morphism
$(X\times_F Y)_L\to X\times_F Y$ is finite,  flat and  radicial: by Lemma
\ref{l4.1a} a) and a limit argument, this implies that the functor is fully faithful. For its essential surjectivity, we steal an idea from \cite[Ch. VIII, \S 1, proof of Th. 2]{lang}. Let $X$ be a 
smooth projective $K$-variety. Then $X$ is defined over a finite sub-extension $L/F$ of
$K/F$. Let
$p^n=[L:F]$, and let
$\Phi_L$ be the absolute  Frobenius of
$L$. The  relative Frobenius morphism (an $L$-morphism)
\[X\to \Phi_L^n X\]
is finite, flat\footnote{To see this, one may use the fact that $X$ is locally isomorphic 
to $\A^n$ for the \'etale topology.} and radicial; by Lemma 
\ref{l4.1a} b), $h(X)\to h(\Phi_L^n X)$ is an
isomorphism in $\Cor_\sim(L,A)$, hence also in $\Cor_\sim(K,A)$. Since
$\Phi_L^n:\Spec L\to \Spec L$ factors through
$\Spec F$, $\Phi_L^n X$ is defined over $F$, proving that the functor is essentially
surjective.
\end{proof}

\subsection{Image motives}\label{s.im-mot} In the study of projective homogeneous varieties, several people (starting with Vishik) have been led to introduce the following 

\begin{defn} Let $X$ be a smooth projective variety. We write
\[\bar \sZ_\sim^*(X,A)=\IM(\sZ_\sim^*(X,A)\to \sZ_\sim^*(X_{F_s},A))\]
where $F_s$ is a separable closure of $F$.
\end{defn}

Using correspondences based on these groups, we define $\overline{\Mot}_\sim(F,A)$, etc. This is mainly interesting when $A=\Z$ or $\Z/p$: for $A=\Q$ the extension of scalars map is injective (transfer argument).

\section{Pure birational motives}\label{s2}

\subsection{First approach: localisation} The first idea to define a notion of pure
birational motives is to localise $\Mot^\eff_\sim(F,A)$ with respect to stable birational morphisms as in \cite{Birat}, hence getting a
functor
\[S_r^{-1}\Sm_\Coprod^\proj(F)\to S_r^{-1}\Mot^\eff_\sim(F,A).\]

\enlargethispage*{30pt}

This idea turns out to be the good one in all important cases, but
to see this we first need some preliminary work. We start by reviewing the sets of morphisms used in \cite[\S 1.7]{Birat}:
\begin{itemize}
\item $S_b$: birational morphisms;
\item $S_h$: projections of the form $X\times (\P^1)^n\to X$;
\item $S_r$ stably birational morphisms: $s\in S_r$ if and only if $s$ is
dominant and gives a purely transcendental function field extension;
\end{itemize}
to which we adjoin
\begin{itemize}
\item $S_b^w$: compositions of blow-ups with smooth centres;
\item $S_r^w=S_b^w\cup S_h$.
\end{itemize}

These morphisms, defined for connected varieties in \cite{Birat}, extend
trivially to the categories of Definition \ref{d1.3} as explained in \cite[Cor. 6.3]{localisation}. More precisely, if $S$ is a set of morphisms of $\Sm(F)$, we define $S^\Coprod\subset\Sm_\Coprod(F)$ as the set of those morphisms which are dominant and whose connected components are all in
$S$. For simplicity, we shall write $S$ rather than $S^\Coprod$ in the sequel. 

By Lemma \ref{l1.3} and \cite[Th. 6.4]{localisation}, the localisation
results of \cite{localisation} and \cite{Birat} extend to the category $\Sm_\Coprod(F)$ and,
moreover, the functors
\[S^{-1}\Sm(F)\to S^{-1}\Sm_\Coprod(F)\]
identify the right hand side with the ``finite coproduct envelope" of the left
hand side. Similarly for their likes with decorations $\Sm^\%$.

We shall view the above morphisms as correspondences via the graph functor. We
introduce two more sets which are convenient here:

\begin{defn}\label{d1.1} We write $\tilde S_b$ and $\tilde S_r$ for the set of
dominant rational maps which induce, respectively, an isomorphism of function
fields and a purely transcendental extension. We let these rational maps act on
pure motives via their graphs, as in \S\ref{s1.4}.
\end{defn}

Thus we have a diagram of inclusions of morphisms on $\Mot_\sim^\eff(F,A)$:
\begin{equation}\label{eq1.1}
\begin{CD}
S_b^w &\quad\subset\quad & S_b^w\cup S_h &\quad=\quad& S_r^w\\
\bigcap & &\bigcap & &\bigcap & \\
S_b &\quad\subset\quad & S_b\cup S_h &\quad\subset\quad & S_r\\
\bigcap & &\bigcap & &\bigcap & \\
\tilde S_b &\quad\subset\quad & \tilde S_b\cup S_h &\quad\subset\quad & \tilde
S_r
\end{CD}
\end{equation}

Let us immediately notice:

\begin{prop}\label{p9.1} Let $S$ be one of the
systems of morphisms in \eqref{eq1.1}. Then the category
$S^{-1}\Mot^\eff_\sim(F,A)$ is an
$A$-linear category provided with a tensor structure, compatible with the
corresponding structures of
$\Mot^\eff_\sim(F,A)$ via the localisation functor.
\end{prop}

\begin{proof} This follows from Theorem \ref{tA.2}, Proposition \ref{pA.2} and the fact that
elements of
$S$ are stable under disjoint unions and products.
\end{proof}

\subsection{Second approach: the Lefschetz ideal}

\begin{defn} We denote by $\sL_\sim$ the ideal of
$\Mot^\eff_\sim(F,A)$ consisting of those morphisms which factor
through some object of the form $P(1)$: this is the \emph{Lefschetz ideal}. It is
a monoidal ideal (\ie, it is closed with respect to composition and tensor products on
the left and on the right).
\end{defn}

\begin{rk} In any additive category $\sA$ there is the notion of
product of two ideals $\sI,\sJ$: 
\[\sI\circ\sJ=\langle f\circ g\mid f\in \sI, g\in \sJ\rangle.\] 

\begin{sloppypar}
If $\sB$ is some given additive subcategory of $\sA$ and
$\sJ=\{f\mid f \text{ factors through some } A\in \sB\}$, then $\sJ$
is idempotent because it is generated by idempotent morphisms,
namely the identity maps of the objects of $\sB$. In
$\sA=\Mot^\eff_\sim(F,A)$, this applies to
$\sL_\sim$.
\end{sloppypar}

On the other hand, in a tensor additive category $\sA$ there is
also the tensor product of two ideals $\sI,\sJ$: for $A,B\in
\sA$
\[(\sI\otimes\sJ)(A,B)=\langle \sA(E\otimes F,B)\circ
\left(\sI(C,E)\otimes \sJ(D,F)\right)\circ
\sA(A,C\otimes D)\rangle\]
where $C,D,E,F$ run through all objects of $\sA$. Coming back to
$\sA=\Mot^\eff_\sim(F,A)$, we have
$\sL_\sim\otimes\sL_\sim=\Mot^\eff_\sim(F,A)(2)\ne \sL_\sim\circ
\sL_\sim=\sL_\sim$. This is in sharp contrast with the case where $\sA$ is rigid
\cite[(6.15)]{ak}.
\end{rk}

\begin{prop}\label{p7.1}  
a) The localisation functor
\[\Mot^\eff_\sim(F,A)\to (S_b^w)^{-1}\Mot^\eff_\sim(F,A)\] 
factors through $\Mot_\sim^\eff(F,A)/\sL_\sim$.\\ 
b)  The functors 
\[\Mot_\sim^\eff(F,A)/\sL_\sim\to
(S_b^w)^{-1}\Mot_\sim^\eff(F,A)\to (S_r^w)^{-1}\Mot_\sim^\eff(F,A)\] 
are both isomorphisms of categories.\\
c) The functor
\[\Mot_\sim^\eff(F,A)/\sL_\sim \to
S_b^{-1}\Mot^\eff_\sim(F,A)\]  
is full.\\
d) For any $s\in \tilde S_r$, $s_*$ becomes invertible in $\tilde
S_b^{-1}\Mot^\eff_\sim(F,A)$.
\end{prop}

\begin{proof} a) By Proposition \ref{p9.1}, it is sufficient to show that
$\bL \mapsto 0$ in 
$(S_b^w)^{-1}\Mot^\eff_\sim(F,A)$. Here as in the proof of b) we shall
use the following formula of Manin \cite[\S 9, Cor. p. 463]{manin}: if
$p:\tilde X\to X$ is a blow-up with smooth centre $Z\subset X$ of
codimension
$n$, then
\begin{equation}\label{eq7.2} h_\sim^\eff(\tilde X)\simeq h_\sim^\eff(X)\oplus
\bigoplus_{i=1}^{n-1} h_\sim^\eff(Z)\otimes \bL^{\otimes i}
\end{equation}
where projecting the right hand side onto $h_\sim^\eff(X)$ we get $p_*$.

In \eqref{eq7.2}, take $X=\P^2$ and for $\tilde X$ the blow-up of
$X$ at (say) $Z=\{(1:0:0)\}$. Since $p$ is invertible in
$(S_b^w)^{-1}\Mot^\eff_\sim(F,A)$, we get $\bL =0$ in this category as
requested.

b) It suffices to show that morphisms of $S_r^w$ become invertible in
$\Mot_\sim^\eff(F,A)/\sL_\sim$, which
immediately follows from \eqref{eq7.2} and the easier projective line
formula.

c) It suffices to show that members of $S_b$ have right inverses in
$\Mot_\sim^\eff(F,A)$: this follows from Lemma
\ref{l7.1}.

d) Let $g:X\ttto Y$ be an element of $\tilde S_r$. Then $X$ is birational
to $Y\times (\P^1)^n$ for some $n\ge 0$, and if $f:X\ttto Y\times
(\P^1)^n$ is the corresponding birational map, its composition
with the first projection $\pi$ is $g$. By Lemma \ref{l1.1}, it suffices to show
that $\pi_*$ is invertible in $\tilde S_b^{-1} \Mot^\eff_\sim(F,A)$, which
follows from b).
\end{proof}

\begin{cor}\label{c7.1} Let $M=\Mot_\sim^\eff(F,A)$.\\
a) The diagram \eqref{eq1.1} 
induces a commutative diagram of categories and functors
\begin{equation}\label{eq1.2}
\begin{CD}
M/\sL_\sim @>\sim>> (S_b^w)^{-1} M @>\sim>>  (S_b^w\cup S_h)^{-1} M @>\sim>>
(S_r^w)^{-1} M\\
&& @V{\text{\rm full}}VV @V{\text{\rm full}}VV @V{}VV \\
&& S_b^{-1} M @>\sim >>  (S_b\cup S_h)^{-1} M @>>>  S_r^{-1} M\\
&& @V{}VV @V{}VV @V{}VV \\
&& \tilde S_b^{-1} M @>\sim>>  (\tilde S_b\cup S_h)^{-1} M @>\sim>>  \tilde
S_r^{-1} M
\end{CD}
\end{equation}
where the functors with a sign $\sim$ are isomorphisms of categories and the
indicated functors are full.\\
*b) If $\car F=0$, all functors are isomorphisms of
categories.
\end{cor}

\begin{proof} a) follows from  Proposition \ref{p7.1}; b) follows from
Hironaka's resolution of singularities (\cf \cite[Lemma
1.7.1]{Birat}).\end{proof}

\begin{rk} Tracking isomorphisms in Diagram \eqref{eq1.2}, one sees that without assuming 
resolution of singularities we get a priori 4 different categories of ``pure
birational motives".  If
$p:\tilde X\to X$ is a birational morphism, then at least $h_\sim(X)$ is a
direct summand of
$h_\sim(\tilde X)$ by Lemma \ref{l7.1}. However it is not clear how to prove
that the other summand is divisible by $\bL $ without using resolution. We shall
get by for special pairs $(A,\sim)$ in Theorem \ref{l5.1} below, using the alteration theorem of de
Jong-Gabber.
\end{rk}

We now introduce:

\begin{defn}\label{d1.2} The category of \emph{pure birational
motives} is
\[\Mot^\b_\sim(F,A)=\left(\Mot^\eff_\sim(F,A)/\sL_\sim\right)^\natural.\]
We also set
\begin{align*}
\Chow^\eff(F,A)&=\Mot_\rat^\eff(F,A)\\
\Chow^\b(F,A)&=\Mot_\rat^\b(F,A).
\end{align*}
When $A=\Z$, we abbreviate this notation to
$\Chow^\eff(F)$ and $\Chow^\b(F)$.\end{defn}

We note:

\begin{prop} Taking pseudo-abelian envelopes, the first functor in Corollary \ref{c7.1} a) induces an isomorphism of categories 
\[\Mot_\sim^\b(F,A)\iso \left((S_b^w)^{-1}\Cor_\sim(F,A)\right)^\natural.\]
 In particular, the functor $(S_b^w)^{-1}\Cor_\sim(F,A)\allowbreak\to (S_b^w)^{-1}\Mot_\sim^\eff(F,A)$ is fully faithful and the functor $\Cor_\sim(F,A)\to S_b^{-1}\Cor_\sim(F,A)$ is full.
\end{prop}

\begin{proof} All follows from Lemma \ref{lA6.1}, except for the last statement which follows from Proposition \ref{p7.1} c).
\end{proof}

In Section \ref{s3}, we shall examine to what extent it is really necessary to
adjoin idempotents in Definition \ref{d1.2}.

\subsection{Third approach: extendible pairs}\label{s2.3} To go 
further, we
need to restrict the adequate equivalence relation we are using:

\begin{defn} \label{d7.1} An adequate pair $(A,\sim)$ is
\emph{extendible} if
\begin{itemize}
\item $\sim$ is defined on cycles over arbitrary quasiprojective
$F$-varieties;
\item it is preserved by inverse image under flat morphisms and
direct image under proper morphisms;
\item if $X$ is smooth projective, $Z$ is a closed subset of $X$
and $U=X-Z$, then the sequence
\begin{equation}\label{eqextend}
\sZ_n^\sim(Z,A)\to \sZ_n^\sim(X,A)\to \sZ_n^\sim(U,A)\to 0
\end{equation}
is exact.
\end{itemize}
\end{defn}

Note that in \eqref{eqextend}, surjectivity always holds because
this is already true on the level of cycles. So the issue is exactness
at $\sZ_n^\sim(X,A)$.

\begin{exs}
a) Rational equivalence (with any coefficients) is extendible.\\
b) Algebraic equivalence (with any coefficients) is extendible, \cf
\cite[Ex. 10.3.4]{fulton}.\\
c) The status of homological equivalence is very interesting:
\begin{enumerate}
\item Under the standard conjecture that
homological and numerical equivalences agree, homological
equivalence with respect to a ``classical" Weil cohomology theory is extendible
if $\car F=0$ (Corti-Hanamura \cite[Prop.
6.7]{ch}). The proof involves resolution of
singularities and the weight spectral sequences for Borel-Moore
Hodge homology, their degeneration at
$E_2$ and the semi-simplicity of numerical motives (Jannsen \cite{jannsen}).
Presumably the same arguments work in characteristic $p$ by using de Jong's
alteration theorem \cite{dJ} instead of Hironaka's resolution of singularities:
we thank Yves Andr\'e for pointing this out.  See \cite[Prop. 1.6]{voisinens} for a more precise statement and a different proof.
\item It seems that the Corti-Hanamura argument implies unconditionally that
Andr\'e's motivated cycles \cite{andre} verify the
axioms of an extendible pair.
\item For Betti cohomology with integral coefficients or $l$-adic cohomology with
$\Z_l$ coefficients, homological equivalence is not extendible. (Counterexample:
$F=\C$, $n=1$,
$Z$ a general surface of degree $\ge 4$ in $\P^3$; this example goes back to Koll\'ar, \cf \cite[p. 134]{kollar}.)  This is closely related to the failure of the Hodge or Tate conjecture integrally for $Z$ (see \cite[\S 2]{soule-voisin}).
\item Hodge cycles with coefficients $\Q$ verify the
axioms of an extendible pair: similarly to (1), the proof involves resolving the singularities of
$Z$ in \eqref{eqextend} and using the semi-simplicity of polarisable pure Hodge
structures. See also Jannsen \cite{jannsen2}.\\ 
We are indebted to Claire Voisin for explaining these last two points.
\item Taking Tate cycles for $l$-adic cohomology, the same argument works if
we assume the semi-simplicity of Galois action on the cohomology of smooth
projective varieties.
\end{enumerate}
\end{exs}

\begin{lemma}\label{l7.4} If $(A,\sim)$ verifies the first two
conditions of Definition \ref{d7.1}, then $(A,\rat)\ge (A,\sim)$
(also over arbitrary quasiprojective varieties).
\end{lemma}

\begin{proof} Again, this follows from \cite[Ex. 1.7.5]{fulton}.\end{proof}

\begin{prop}\label{l4.2'} Let $(A,\sim)$ be an extendible pair. For
two smooth projective varieties $X,Y$, let $\sI_\sim(X,Y)$ be
the subgroup of $\sZ^{\dim Y}_\sim(X\times Y,A)$ consisting of those
classes vanishing in $\sZ^{\dim Y}_\sim(U\times Y,A)$ for some open
subset
$U$ of $X$. Then $\sI_\sim$ is a monoidal ideal in $\Cor_\sim(F,A)$. 
\end{prop}

\begin{proof} Note that by Lemma \ref{l7.4} and the third condition of
Definition \ref{d7.1}, the map $\sI_\rat(X,Y)\to \sI_\sim(X,Y)$
is surjective for any $X,Y$: this reduces us to the case
$\sim=\rat$. We further reduce immediately to $A=\Z$.

Let
$X,Y,Z$ be 3 smooth projective varieties. If
$U$ is an open subset of $X$, it is clear that the usual formula
defines a composition of correspondences 
\[CH^{\dim Y}(U\times Y)\times CH^{\dim Z}(Y\times
Z)\to CH^{\dim Z}(U\times Z)\] 
and that this composition commutes with restriction to smaller and
smaller open subsets. Passing to the limit on
$U$, we get a composition 
\[CH^{\dim Y} (Y_{F(X)})\times CH^{\dim
Z}(Y\times Z)\to CH^{\dim Z}(Z_{F(X)})\] 
or 
\[CH_0(Y_{F(X)})\times CH^{\dim Z}(Y\times Z)\to
CH_0(Z_{F(X)}).\]

Here we used the fact that (codimensional) Chow groups commute with
filtering inverse limits of schemes, see \cite{bloch}.

We now need to prove that this pairing factors through 
$CH_0(Y_{F(X)})\times CH^{\dim Z}(V\times Z)$ for any open subset $V$ of $Y$.
One checks that it is induced by the standard action of
correspondences in $CH^{\dim Z}(Y_{F(X)}\times_{F(X)} Z_{F(X)})$ on groups of
$0$-cycles. Hence it is sufficient to show that the standard action of
correspondences factors as indicated, and up to changing the base field we may
replace $F(X)$ by $F$.

We now show that the pairing
\[CH_0(Y)\times CH^{\dim Z}(Y\times Z)\to
CH_0(Z)\]
factors as indicated. The proof is a variant of Fulton's proof of
the Colliot-Th\'el\`ene--Coray theorem that $CH_0$ is a birational invariant
of smooth projective varieties \cite{ctc}, \cite[Ex. 16.1.11]{fulton}.
Let $M$ be a proper closed subset of $Y$, and $i:M\to Y$ be the
corresponding closed immersion. We have to prove that for any
$\alpha\in CH_0(Y)$ and 
$\beta\in CH_{\dim Y}(M\times Z)$, 
\[(i\times 1_Z)_*(\beta)(\alpha):= (p_2)_*((i\times 1_Z)_*\beta\cdot 
p_1^*\alpha)=0\]
where $p_1$ and $p_2$ are respectively the first and second
projections on $Y\times Z$.

We shall actually prove that $(i\times 1_Z)_*\beta\cdot
p_1^*\alpha=0$. For this, we may assume that $\alpha$ is
represented by a closed point $y\in Y$ and $\beta$ by some integral
variety $W\subseteq M\times Z$. Then $(i\times 1_Z)_*\beta\cdot
p_1^*\alpha$ has support in $(i\times 1_Z)(W)\cap (\{y\}\times Z)
\subset (M\times Z)\cap (\{y\}\times Z)$. If $y\notin M$, this
subset is empty and we are done. Otherwise, up to rational
equivalence, we may replace $y$ by a $0$-cycle
disjoint from $M$ (\cf \cite{roberts}), and we are back to the
previous case.

This shows that $\sI_\sim$ is an ideal of $\Cor_\sim(F,A)$. The fact
that it is a monoidal ideal is essentially obvious.\end{proof}

\begin{defn}\label{d1.2a} For an extendible pair $(A,\sim)$, we abbreviate the
notation $\Cor_\sim(F,A)/\sI_\sim$ (resp.  $\left(\Mot^\eff_\sim(F,A)/\sI_\sim\right)^\natural$) into $\Cor^\o_\sim(F,A)$ (resp. $\Mot^\o_\sim(F,A)$). ($\o$ stands for ``open''.) We write $h^\o_\sim(X)$ for the image of $h_\sim(X)$ in $\Mot^\o_\sim(F,A)$. We also set
$\Chow^\o(F,A)=\Mot_\rat^\o(F,A)$ and $\Chow^\o(F)=\Chow^\o(F,\Z)$.
\end{defn}

For future reference, let us record here the value of the Hom groups
in the most important case, that of rational equivalence (see also
Remark \ref{rhomsim} 2) below):
\begin{lemma}\label{ehomrat} We have
\[\Cor_\rat^\o(F,A)([X],[Y])= CH_0(Y_{F(X)})\otimes A.\]
\end{lemma}

\begin{prop}\label{l7.3} In $\Cor^\o_\sim(F,A)$,\\
a) $(g\circ f)_*=g_*\circ f_*$ for any composable rational maps
$X\stackrel{f}{\ttto}Y\stackrel{g}{\ttto}Z$.\\
b) \cite[Ex. 16.1.11]{fulton}  $f^*f_*=1_X$ and $f_*f^*=1_Y$  for any
birational map $f:X\ttto Y$.\\ 
c) Morphisms of $\tilde S_r$ (see Definition \ref{d1.1}) are invertible.
\end{prop}

\begin{proof} a) Let $F$ be the fundamental set of $f$, $G$ be the
fundamental set of $g$, $U=X-F$,
$V=Y-G$. By assumption, $f(U)\cap V\neq\emptyset$, hence
$W=f^{-1}(V)$ is a nonempty open subset of $U$, on which $g\circ f$
is a morphism.

Let us abuse notation and still write $f$ for the morphism $f_U$,
etc. Then, by definition
\[g_*\circ f_*=(p_{XZ})_*((\bar\Gamma_f\times Z)\cap (X\times
\bar\Gamma_g))\] 
(note that the two intersected cycles are in good
position). This cycle clearly contains $(g\circ f)_*=\bar\Gamma_{g\circ
f}$ as a closed subset. One sees immediately that the restriction of
$g_*\circ f_*$ and $(g\circ f)_*$ to $W\times Z$ are equal.

b) is proven in the same way (or is a special case of a)).

c) Let $g:X\ttto Y$ be an element of $\tilde S_r$. Then $X$ is birational
to $Y\times (\P^1)^n$ for some $n\ge 0$, and if $f:X\ttto Y\times
(\P^1)^n$ is a birational map, its composition
with the first projection $\pi$ is $g$. By a) and b), it suffices to show
that $\pi_*$ is invertible in $\Cor_\sim(F,A)/\sI_\sim$. For this we may
reduce to $n=1$ and even to $Y=\Spec F$ since $\sI_\sim$ is a monoidal
ideal. Let $s:\Spec F\to \P^1$ be the $\infty$ section: it suffices to
show that $(s\circ \pi)_*=1_{\P^1}$. But the cycle 
$(s\circ\pi)_*-1_{\P^1}$ on $\P^1\times\P^1$ is linearly equivalent to
$\infty\times\P^1$ (this is the idempotent defining the Lefschetz motive),
and the latter cycle vanishes when restricted to $\A^1\times\P^1$.
 \end{proof}

We shall also need the following lemma in the proof of Proposition 
\ref{ph0} c).

\begin{lemma}\label{l2.3} Let $L/K$ be an extension of function fields 
over $F$, with 
$K=F(X)$ and $L=F(Y)$ for $X,Y$ two smooth projective $F$-varieties. Let 
$\phi:Y\ttto X$ be the rational map corresponding to the inclusion $K\inj 
L$. Let $Z$ be another smooth projective $F$-variety. 
Then the map
\[\Chow^\o(F,A)(h^\o(X),h^\o(Z))\to \Chow^\o(F,A)(h^\o(Y),h^\o(Z))\]
given by composition with $\phi_*:h^\o(Y)\to h^\o(X)$ (see \ref{s1.4}) 
coincides via Lemma \ref{ehomrat} with the base-change map $CH_0(Z_K)\otimes 
A\to CH_0(Z_L)\otimes A$.
\end{lemma}

\begin{proof} Let $V\subseteq Y$ and $U\subseteq X$ be open subsets such 
that $\phi$ is defined on $V$ and $\phi(V)\subseteq U$. Up to shrinking $U$, we 
may assume that $\phi$ is flat \cite[11.1.1 ]{EGAIV}. As in the proof of 
Proposition \ref{l4.2'}, the composition of correspondences induces a 
pairing
\[CH^{\dim X}(V\times U)\times CH^{\dim Z}(U\times Z)\to CH^{\dim 
Z}(V\times Z)\]
and the action of $\phi_*\in CH^{\dim X}(V\times U)$ on $\alpha\in 
CH^{\dim 
Z}(U\times Z)$ is given by the flat pull-back of cycles. Therefore, 
$\phi_*$ 
induces in the limit the flat pull-back of $0$-cycles from $CH_0(Z_K)$ to 
$CH_0(Z_L)$.
\end{proof}

\begin{rks}\label{rhomsim} 1) Propositions \ref{l4.2'} and \ref{l7.3} a)
  were independently 
observed by Markus Rost in the case $\sim=\rat$ \cite[Prop. 3.1 and Lemma
3.3]{merk}. We are indebted to Karpenko for pointing this out and for
referring us to Merkurjev's preprint \cite{merk}.

2) In $\Cor^\o_\sim(F,A)$, morphisms are by definition given by the
formula
\[\Cor^\o_\sim(F,A)([X],[Y])=\colim_{U\subseteq
X} \sZ^{\dim Y}_\sim(U\times Y,A).\]

The latter group maps onto $\sZ_0^\sim(Y_{F(X)},A)$. If
$\sim=\rat$, this map is an isomorphism (see Lemma \ref{ehomrat}). For
other equivalence relations, this
is far from being the case: for example, if $\sim=\alg$, $F$ is
algebraically closed,
$X,Y$ are two curves and (say) $A=\Z$, then
\begin{multline*}
\sZ^1_\alg(X\times Y,\Z)=NS(X\times Y)=NS(X)\oplus
NS(Y)\oplus
\Hom(J_X,J_Y)\\=\Z\oplus\Z\oplus \Hom(J_X,J_Y)
\end{multline*}
where $NS$ is the N\'eron-Severi group and $J_X,J_Y$ are the Jacobians of
$X$ and $Y$. On the other hand,
\[\sZ_0^\alg(Y_{F(X)},\Z)= NS(Y_{F(X)})=\Z.\]

When we remove a point from $X$, we kill the factor $NS(X)=\Z$. But any
two points of $X$ are algebraically equivalent, so removing further
points does not modify the group any further. Hence
\[\colim_{U\subseteq X} \sZ^{\dim Y}_\alg(U\times Y,\Z)
=\Z\oplus\Hom(J_X,J_Y).\]

We thank Colliot-Th\'el\`ene for helping clarify this matter.
\end{rks}

\subsection{The main theorem}
We now extend the ideal $\sI_\sim$ from $\Cor_\sim(F,A)$
to $\Mot^\eff_\sim(F,A)$ in the usual way (\cf \cite[Lemme 1.3.10]{ak}),
without changing notation. By Propositions \ref{p7.1} a) and 
\ref{l7.3}, we get a composite functor
\begin{equation}\label{eq7.4}
\Mot_\sim^\b(F,A)\to (\tilde S_r^{-1}\Mot_\sim^\eff(F,A))^\natural \to 
\Mot_\sim^\o(F,A)
\end{equation}
for any extendible pair $(A,\sim)$. Since both categories are (idempotent completions of) full images
of $\Mot_\sim^\eff(F,A)$, this functor is automatically full. We are
going to show that it is an equivalence of categories in some important
cases.

\begin{thm}\label{l5.1}  Let $(A,\sim)$ be an extendible pair. Suppose that
the exponential characteristic $p$ of $F$ is invertible in $A$. Then the functor \eqref{eq7.4} is an isomorphism of categories.
\end{thm}

\begin{proof}\footnote{We thank N. Fakhruddin for his help, which removes the 
recourse to Chow's moving lemma in \cite{birat}.} We have to show that
$\sI_\sim(M,N)\subseteq
\sL_\sim(M,N)$ for any $M,N\in \Mot_\sim^\eff(F,A)$. Proposition \ref{l3.1a} reduces us to the case where $F$ is \emph{perfect}. 
Clearly we may assume $M=h_\sim(X)$, $N=h_\sim(Y))$ for two smooth
projective varieties $X,Y$. 

Let $f\in \sI_\sim(h_\sim(X),h_\sim(Y))$. By the third condition in
Definition \ref{d7.1}, the cycle class
$f\in\sZ_{\dim X}^\sim(X\times Y,A)$ is of the form
$(i\times 1_Y)_*g$  for some closed immersion $i:Z\to X$, where $g\in
\sZ_{\dim X}^\sim(Z\times Y,A)$. Let $\tilde g$ be a cycle representing $g$.
Write $\tilde g = \sum_k a_k g_k$, with $a_k\in A$ and $g_k$ irreducible. Then
$(i\times 1_Y)_*(g_k)\in\sI_\sim(h_\sim(X),h_\sim(Y))$. This reduces us to the
case where $g$ is represented by an irreducible cycle $\tilde g$.

Choose $Z$ minimal among the closed subsets of $X$ such that $\tilde g$ is
supported on $Z\times Y$. In particular, $Z$ is irreducible. 

Consider $Z$ with its reduced structure. Let $l$ be a prime number different from $p$: by Gabber's refinement of de Jong's theorem \cite[Th. X 2.1]{gabber}, we may choose a
proper, generically finite morphism $\pi_l:\tilde Z_l\to Z$  where
$\tilde Z_l$ is smooth projective (irreducible) and
$\pi_l$ is an alteration of generic degree $d_l$ prime to $l$. (Recall that an alteration is a proper, generically finite morphism.)

By the minimality of $Z$, the support of $\tilde g$ has nonempty intersection
$\tilde g_1$ with $V\times Y$, where $V = Z-(Z_{sing}\cup T)$ with $Z_{sing}$
the singular locus of $Z$ and $T$ the closed subset over which $\pi_l$ is not
finite. Let $\pi_V:\pi_l^{-1} (V)\to V$ be the map induced by $\pi_l$:  note that $\pi_V$ is flat since $V$ and $\pi_l^{-1}(V)$ are smooth. We then have an equality of cycles 
\[d_l\tilde g_1 =(\pi_V\times 1_Y)_*(\pi_V\times 1_Y)^* \tilde g_1.\]

Let $\gamma_l$ be the closure of $(\pi_V\times 1_Y)^* \tilde g_1$ in $\tilde Z_l$\footnote{More correctly, the cycle associated to the schematic closure of $(\pi_V\times 1_Y)^{-1} (\tilde g_1)$ in $\tilde Z_l$: take the topological closure of each component of $(\pi_V\times 1_Y)^* \tilde g_1$ and keep the same multiplicities.}. We get an equality of cycles (the support of $(\pi_V\times 1_Y)_*(\pi_V\times 1_Y)^* \tilde g_1$ is dense in that of $(\pi_l\times 1_Y)_*\gamma_l$): 
\[d_l\tilde g =(\pi_l\times 1_Y)_*\gamma_l.\]

Let $d=gcd_l(d_l)$, which is a power of $p$; then $d=gcd(d_{l_1},\dots, d_{l_r})$ for some finite set of primes $\{l_1,\dots,l_r\}$. For simplicity, write $Z_{l_i}=Z_i$, $\pi_{l_i}=\pi_i$  and $\gamma_{l_i}=\gamma_i$. 

Let $h_i=d^{-1} [\gamma_i]\in \sZ_{\dim
X}^\sim(\tilde Z_i\times Y,A)$. Choose $a_1,\dots, a_r\in \Z$ such that $d=\sum_i a_i d_i$, so that 
\[f=\sum_i a_i ((i\circ \pi_i)\times 1_Y)_*h_i. \]

Then the correspondence
$f\in
\Mot^\eff_\sim(F)(h_\sim(X),h_\sim(Y))$ factors as
\[\begin{CD}
h_\sim(X)@>(i\circ \pi)^*>> h_\sim(\coprod\tilde Z_i)(\dim X -\dim
Z)@>(h_i)>> h_\sim(Y)
\end{CD}\]
(see \eqref{eq7.3}), which concludes the proof.
\end{proof}

\begin{cor}\label{c7.2} Under the assumptions of Theorem \ref{l5.1},
all the categories of Diagram \eqref{eq1.2} are isomorphic to
$\Mot_\sim^\eff(F,A)/\sI_\sim$.
\end{cor}

\begin{proof} By Proposition \ref{p7.1} b) and d) we already know that the
categories
$\Mot_\sim^\eff(F,A)/\sL_\sim$,
$(S_b^w)^{-1}\Mot_\sim^\eff(F,A)$ and 
$(S_r^w)^{-1}\Mot_\sim^\eff(F,A)$ are isomorphic and that $(\tilde
S_b)^{-1}\Mot_\sim^\eff(F,A)$ and $(\tilde
S_r)^{-1}\Mot_\sim^\eff(F,A)$ are isomorphic. We also know that the
functor
$\Mot_\sim^\eff(F,A)/\sL_\sim\to
(S_b)^{-1}\Mot_\sim^\eff(F,A)$ is full (Proposition \ref{p7.1} c)): by Theorem \ref{l5.1}, this implies that it is an isomorphism. To conclude the proof, it is
sufficient to show that any morphism of
$\tilde S_r$, hence of $S_r$, has a right inverse in
$\Mot_\sim^\eff(F,A)/\sL_\sim$ (see \eqref{eq1.2}). Since $\tilde S_r$ is generated by $\tilde
S_b$ and projections of the form $X\times
\P^1\to X$ (\cf proof of Proposition \ref{p7.1} d)) and since this
is obvious for these projections, we are left to prove it for
elements $f:X\ttto Y$ of $\tilde S_b$. But we have $f_*f^*=1_X$ in
$\Mot_\sim^\eff(F,A)/\sI_\sim$  by Proposition \ref{l7.3} b), hence in
$\Mot_\sim^\eff(F,A)/\sL_\sim$ by Theorem \ref{l5.1}. \end{proof}

\subsection{Birational image motives} \label{s.im-mot-bir} Based on the categories of Subsection \ref{s.im-mot}, we define categories $\overline{\Mot}_\sim^\b(F,A)$. If $\sim$ is extendible and $p$ is invertible in $A$, the analogue of Theorem \ref{l5.1} holds, with the same proof.

\subsection{Recapitulation, comments and notation} In Definition \ref{d1.2}, we associated to any admissible pair $(A,\sim)$ a category of birational motives $\Mot^\b_\sim(F,A)$. If $(A,\sim)$ is extendible (Definition \ref{d7.1}), we introduced in Definition \ref{d1.2a}  another category $\Mot^\o_\sim(F,A)$ plus a full functor $\Mot^\b_\sim(F,A)\to \Mot^\o_\sim(F,A)$. We showed in Theorem \ref{l5.1} that this functor is an isomorphism of categories when the exponential characteristic $p$ is invertible in $A$; in particular, this is true for any $A$ in characteristic $0$. This gives a great flexibility in computing Hom groups, as in some cases one can use their ``algebraic'' description in terms of killing the Lefschetz motive, and in other cases their ``geometric'' description as Chow groups of $0$-cycles if $\sim$ is rational equivalence. 

In the sequel, we commit the abuse of notation which consists of writing $\Mot_\sim^\o$ for $\Mot_\sim^\b$ even when we don't know if the pair \text{$(A,\sim)$} is extendible (notably, when $\sim$ is numerical equivalence). We do this because we feel that keeping the distinction would create more confusion than this choice.

\section{Examples}\label{sexamples}

We give some examples and computations of birational motives.

\subsection{Varieties with trivial birational motive}\label{s3.1} They were initially studied by Bloch-Srinivas \cite{bs} over a universal domain.  The reader should compare the following to \cite[Th. 8.5.1]{Birat}.

\begin{prop}\label{p3.4} Let $A$ be a connected commutative ring, and let $X$ be a smooth projective $F$-variety. Then the following conditions are equivalent:
\begin{thlist}
\item For any smooth projective  $F$-variety $Y$, $CH_0(X_{F(Y)})\otimes A\iso A$ (by the degree map).
\item $CH_0(X_{F(X)})\otimes A\iso A$. 
\item The class of the generic point $\eta_X$ in $CH_0(X_{F(X)})\otimes A$ belongs to $\IM(CH_0(X)\otimes A\to CH_0(X_{F(X)})\otimes A)$.
\item $h^\o(X)=\un$ in $\Chow^\o(F,A)$.
\item (For $A=\Z$:) $M_0(F)\iso A^0(X,M_0)$ for any cycle module $M$.
\end{thlist}
If $p$ is invertible in $A$, they are also equivalent to
\begin{thlist}
\item[\rm (vi)] For any extension $K/F$, $CH_0(X_K)\otimes A\iso A$.
\end{thlist}
If $F$ is a universal domain and $A\supseteq \Q$, they are also equivalent to
\begin{thlist}
\item[\rm (vii)] $CH_0(X)\otimes A\iso A$.
\item[\rm (viii)] $CH_0(X)\iso \Z$.
\end{thlist}
\end{prop}

(Parts of this proposition are standard, see \eg \cite[Lemma 1.3]{auel}.)

\begin{proof} (i) $\Rightarrow$ (ii) $\Rightarrow$ (iii) are obvious. 
By Lemma \ref{ehomrat}, the map of (iii) can be translated into 
\[\Chow^\o(F,A)(\un,h^\o(X))\to \Chow^\o(F,A)(h^\o(X),h^\o(X))\]
via the projection $h^\o(X)\to h^\o(\Spec k)=\un$. Since $\eta_X$ represents the identity endomorphism of $h^\o(X)$, (iii) means that the latter factors through $\un$. Since $\End(\un)=A$, the resulting idempotent endomorphism of $\un$ must be $0$ or $1$; so $h^\o(X)=0$ or $\un$, but the first case is impossible as it would imply that $\eta_X=0$, while $\deg(\eta_X)=1$. So (iii) $\Rightarrow $ (iv). Using Lemma \ref{ehomrat} again, we get (iv) $\Rightarrow$ (i). 

(vi) $\Rightarrow$ (i) is obvious; to prove the converse, we reduce to $F$ perfect by using Proposition \ref{l3.1a}, and then to $K/F$ finitely generated by a limit argument. Then $K$ is the function field of some smooth $F$-variety. We argue as in the proof of Theorem \ref{l5.1}: using  \cite[Th. X.2.1]{gabber}, we can find finite extensions $L_i/K$ such that $L_i=F(Y_i)$ for $Y_i$ smooth projective, such that the gcd of the $[L_i:K]$'s is a power of $p$. Then $(CH_0(X_K)\otimes A)_{\deg =0}$ is a direct summand of $\bigoplus_i (CH_0(X_{L_i})\otimes A)_{\deg =0}=0$ by a transfer argument, hence (vi).

(iv) $\Rightarrow$ (v) $\Rightarrow$ (iii): see Section \ref{s.6}.

It remains to prove (iii) $\Leftarrow$ (vii) $\Rightarrow$ (viii) when $F$ is a universal domain, as (viii) $\Rightarrow$ (vii) is obvious.. The implication (vii) $\Rightarrow$ (iii) is the classical Bloch-Srinivas argument \cite[Prop. 1]{bs}: $X$ is defined over a subfield $F'\subset F$ finitely generated over the prime field; for clarity, write $X'$ for this $F'$-model. Now $F'(X')$ embeds into $F$ over $F'$. Since $\Ker(CH_0(X'_{F'(X')})\to CH_0(X'_F)=CH_0(X))$ is torsion by a transfer argument, (vii) implies that $CH_0(X'_{F'(X')})\otimes A\iso A$. Thus $\eta_{X'}$ is $A$-rationally equivalent to a closed point of $X'$, hence (iii). If (vii) is true, then $\Albv(X)(F)\otimes A=0$ where $\Albv(X)$ is the Albanese variety of $X$; this implies $\Albv(X)=0$. But Ro\v\i tman's theorem \cite{roitman2} then implies that $CH_0(X)_\tors=0$, whence (viii).
\end{proof}

\begin{cor}\label{c3.3} Conditions (i)--(v) of Proposition \ref{p3.4} are stable under products of varieties; so are (vi), (vii) and (viii) under the stated conditions on $A$ and $F$.
\end{cor}

\begin{proof} Indeed, this is obviously the case for Condition (iv).
\end{proof}

\begin{rks} 
1) Condition (v) of Proposition \ref{p3.4} can be extended to any $A$ if we consider cycle modules with coefficients in $A$.\\
2) Except for (iv), Corollary \ref{c3.3} can also be proven without reference to birational motives when $A\supseteq \Q$, using that the product map
\[(CH_0(X)\otimes A)\otimes (CH_0(Y)\otimes A)\to CH_0(X\times Y)\otimes A\]
is then surjective for any smooth projective $X,Y$: reduce to $F$ algebraically closed by a transfer argument, when this even holds integrally.
\end{rks}

We now give some examples. In Part 3 of the following proposition, the Betti numbers $b^i(X)=\dim H^i(X)$ refer to a ``classical'' Weil cohomology $H$: Betti or de Rham in characteristic $0$, crystalline in characteristic $>0$, $l$-adic in characteristic $\ne l$. It is known that $b^i(X)$ does not depend on the choice of such a Weil cohomology.

\begin{prop}\label{trcc}  1) If $X$ is retract rational, then $h^\o(X)=\un$ in $\Chow^\o(F,\Z)$.\\
2) If $X$ is rationally chain connected, then $h^\o(X)=\un$ in $\Chow^\o(F,\Q)$.\\
3) If $h^\o(X)=\un$ in $\Chow^\o(F,\Q)$, then $b^1(X)=0$ and $b^2(X)=\rho(X)$ (Picard number). \\
4) If $\dim X=2$, the converse of 3) is true if and only if $X$ verifies Bloch's conjecture on $0$-cycles. 
\end{prop}

\begin{proof} 1) This follows from \cite[Prop. 8.6.2]{Birat} and the functor \eqref{eq5.1} below. (One could also give a direct proof.)

2) Let $\overline{F(X)}$ be an algebraic closure of $F(X)$: then $X(\overline{F(X)})/R=*$. Since the group of $0$-cycles on $X_{\overline{F(X)}}$ is
generated by $X(\overline{F(X)})$, this in turn implies that $CH_0(X_{\overline{F(X)}})\iso \Z$,
which implies by a transfer argument that $CH_0(X_{F(X)})\otimes\Q\iso\Q$.

3) Since the hypothesis and conclusion do not change by extension of $F$, we may assume that $F$ is a universal domain. We  use Theorem \ref{l5.1}: in $\Chow^\eff=\Chow^\eff(F,\Q)$ we get a decomposition
\[h(X)=\un\oplus M\otimes \bL\]
for some $M\in \Chow^\eff$. Applying the cycle class map, we get a commutative diagram
\[\begin{CD}
CH^1(X)\otimes K@= CH^0(M)\otimes K\\
@V\cl_X^1VV @V\cl_M^0 VV\\
H^2(X)@= H^0(M). 
\end{CD}\]

Here $K$ is the field of coefficients of $H$ and, as usual, $CH^i(M):= \Chow^\eff(M,\bL^i)$ (giving back the rational Chow groups of smooth projective varieties) and $\cl$ is the cycle class map; for simplicity, we neglect Tate twists on cohomology. But $\cl_M^0$ is an isomorphism, as one sees by writing $M$ as a direct summand of $h(Y)$ for some smooth projective $Y$; therefore $\cl^1_X$ is an isomorphism as well. Since this map factors through the N\'eron-Severi group $\NS(X)\otimes K$, this implies $\Pic^0(X)=0$ (hence $b^1(X)=0$), and $b^2(X)=\rho(X)$ as requested.

4) The conditions in the conclusion of 3) imply $\Albv(X)=0$ and (under Bloch's conjecture) $T(X_K)=0$ for any extension $K/F$, where $T$ is the Albanese kernel; the conclusion now follows from Condition (i) of Proposition \ref{p3.4}.
\end{proof}

\begin{rks}\label{r3.1} 1) As noted in \cite[Ex. 7.3]{FqX}, an Enriques surface verifies the conditions of Proposition \ref{p3.4} (for $2$ invertible in $A$); this can be recovered from Proposition \ref{trcc} 4) in a rather silly way. On the other hand, Inose-Mizukami and Voisin's proofs of the Bloch conjecture for some quotients of hypersurfaces by finite groups \cite{in-miz,godeaux} give examples of surfaces of general type having trivial birational motive (with $\Q$-coefficients), which shows once again how motivic information is in some sense orthogonal to geometric information related to the Kodaira dimension.  For a more refined example, see Remark 3.

2) Applying the reasoning in the proof of Proposition \ref{trcc} 3) to $CH^2$ and $CH_1$, one recovers some of the representability results of \cite{bs} in a different way. (The situation considered by Bloch and Srinivas is more general, and in the present terms amounts to the following: assume that, in $\Chow^\o(F,\Q)$, $h^\o(X)$ is isomorphic to a direct summand of $h^\o(Y)$ for some smooth projective variety $Y$ of dimension $n\le 3$.)

3) Let $X$ be a smooth projective variety such that $h^\o(X)=\un$ in $\Chow^\o(F,\Q)$. For simplicity, assume that $X$ has a rational point $x$. By Condition (iii) of Proposition \ref{p3.4}, there is an integer $N>0$ such that $N(\eta_X-x)=0$ in $CH_0(X_{F(X)})$. Then in $\Chow^\o(F,\Z)$, we have
\[h^\o(X)=\un \oplus M \qquad \text{with } N 1_M=0.\]

Indeed, $x$ defines an idenpotent endomorphism of $h^\o(X)$ which splits off the summand $\un$, and $\eta_X-x$ is the complementary idempotent. It follows that $N CH_0(X_K)_0=0$ for any extension $K/F$ and (for instance) that 
\[N\Coker(M_n(K)\to A^0(X_K,M_n))=N\Ker(A_0(X_K,M_n)\to M_n(K))=0\] 
for any cycle module $M$ and any $K\supseteq F$ (see \S \ref{s.6}): compare \cite[Th. 1.4]{auel}. 

 If  $N$ is minimal, then $N>1$ is an obstruction to having $h^\o(X)=\un$ in $\Chow^\o(F,\Z)$: this obstruction has been studied recently in \cite{auel}, \cite{voisin2} and \cite{voisin3}. Using the cycle module $M_n(K)=H^{n}(K,\Q/\Z(n-1))$ for $n=1$, one finds that $N$ is divisible by the exponent $e$ of $H^1_\et(X_{\bar F},\Q/\Z)$.  One can show that  $N=e$ if $F$ is algebraically closed and $X$ is a surface \cite{exponent}; for $e=1$, this was proven by Voisin in \cite[Prop. 2.2]{voisin2} and by Auel, Colliot-Th\'el\`ene and Parimala in \cite[Cor. 1.10]{auel}. For example, $N=2$ for an Enriques surface and $N=1$ for Barlow's surface (of general type) \cite{barlow1,barlow2}, showing that its motive is $\un$ in $\Chow^\o(F,\Z)$. (See the recent survey paper \cite{pg0survey} for more examples of surfaces of general type with $p_g=0$.)
\end{rks}

\subsection{Quadrics} Suppose $\car F\ne 2$ and let $X$ be a smooth projective quadric over
$F$. By a theorem of Swan and Karpenko \cite{sw,kar}, the degree map
\[\deg:CH_0(X)\to \Z\]
is injective, with image $\Z$ if $X$ has a rational point and $2\Z$ otherwise. This implies:

\begin{prop}\label{p3.8} Let $X,Y$ be two smooth projective over $F$. Suppose that $Y$ is a quadric. Then, in
$\Chow^\o(F)$, we have
\[\Hom( h^\o(X), h^\o(Y))=
\begin{cases}
\Z &\text{if $Y_{F(X)}$ is isotropic}\\
2\Z&\text{otherwise}
\end{cases}
\]
where we have used the degree map $\deg:CH_0(Y_{F(X)})\to\Z$. Similarly, in
$\overline{\Chow}^\o(F,\Z/2)$ (see \S \ref{s.im-mot-bir}), we have
\[\Hom( h^\o(X), h^\o(Y))=
\begin{cases}
\Z/2 &\text{if $Y_{F(X)}$ is isotropic}\\
0&\text{otherwise.}
\end{cases}
\]
\end{prop}

\begin{rk} Much work has been done recently on torsion in $CH_0$ of projective homogeneous varieties: we may quote  \cite{cgm,krashen,psz,cm}\dots\ There are many examples of projective homogeneous varieties other than quadrics for which $CH_0(Y)$ is torsion-free; by \cite[Cor. 4.3]{cm}, this is always the case if $Y$ is isotropic. This allows one to extend the second part of Proposition \ref{p3.8} to arbitrary projective homogeneous $Y$'s (with suitable coefficients). On the other hand, there are examples of anisotropic $Y$'s such that $CH_0(Y)_\tors\ne 0$ (\cite[Prop. 1.1]{krashen}, \cite[\S 18]{cm}), so the first part of Proposition \ref{p3.8} does not extend in full generality.
\end{rk}

\subsection{The nilpotence conjecture} It is:

\begin{conj}\label{c3.1} For any two adequate pairs $(A,\sim),(A,\sim')$ with $A\supseteq \Q$
and $\sim\ge \sim'$, and any $M\in \Mot_\sim(F,A)$, $\Ker(\End(M)\to \End(M_{\sim'}))$ is
nilpotent. (We say that the kernel of $\Mot_\sim(F,A)\allowbreak\to
\Mot_{\sim'}(F,A)$ is \emph{locally nilpotent}.)
\end{conj}

Since $\rat$ is the finest (\resp $\num$ is the coarsest) adequate equivalence relation, this
conjecture is clearly equivalent to the same statement for $\sim =\rat$ and $\sim'=\num$, but
it may be convenient to consider it for selected adequate equivalence relations. For example:

\begin{prop}\label{p3.2} a) Conjecture \ref{c3.1} is true for $M\in \Mot_\sim^\eff(F,A)$ (and any $\sim'\le \sim$) provided $M$ is
finite-dimensional in the sense of Ki\-mu\-ra-O'Sullivan \cite[Def. 3.7]{kimura}. In
particular, it is true if $M$ is of abelian type, \ie $M$ is a direct summand of $h_\sim(A_K)$ for
$A$ an abelian $F$-variety and $K$ a finite extension of $F$.\\ 
b) If $\sim=\hom$,
$\sim'=\num$, the condition of a) is equivalent to the
\begin{quote}
\emph{Sign conjecture}: if $H$ is the Weil cohomology theory defining $\hom$, the projector of
$\End H(M)$ projecting $H(M)\allowbreak=H^+(M)\oplus H^-(M)$ onto its summand $H^+(M)$ is algebraic.
\end{quote}
 In
particular, it is true if $M$ satisfies the Standard conjecture C (algebraicity of the
K\"unneth projectors).\\ c) Conjecture \ref{c3.1} is true in the following cases:
\begin{thlist}
\item $\sim=\rat$, $\sim'=\tnil$;
\item $\sim=\rat$, $\sim'=\alg$.
\end{thlist}
\end{prop}

\begin{proof} a) This is a theorem of Kimura and O'Sullivan, \cf  \cite[Prop. 7.5]{kimura}, 
\cite[Prop. 9.1.14]{ak}. The second assertion follows from Kimura's results, \cf
\cite[Ex. 7.6.3 4)]{kmp}. b) See
\cite[Th. 9.2.1 c)]{ak}. c) (i) follows from the Voevodsky-Kimura lemma that smash-nil\-pot\-ent
correspondences are nilpotent,
\cf
\cite[Lemma 2.7]{voenil}, \cite[Prop. 2.16]{kimura},
\cite[Lemma 7.4.2 ii)]{ak}. (ii) follows from (i) and Voevodsky's theorem that $\alg\ge \tnil$,
\cite[Cor. 3.2]{voenil}.
\end{proof}

Let us recall some conjectures which imply Conjecture \ref{c3.1}:

\begin{prop} \label{p3.3}
a) Conjecture \ref{c3.1} is implied by Voevodsky's conjecture that smash-nilpotence
equivalence equals numerical equivalence \cite[Conj. 4.2]{voenil}.\\ 
b) It is also implied by the sign conjecture plus the Bloch-Bei\-lin\-son--Murre conjecture \cite{jannsen2,murre2}.
\end{prop}

\begin{proof}   a)  This follows from Proposition \ref{p3.2} c) (i). 
b)  Recall that the
Bloch-Beilinson conjecture is equivalent to Murre's conjecture in 
\cite{murre2} by \cite[Th.
5.2]{jannsen2}. Now the formulation of the former conjecture, 
\cite[Conj. 2.1]{jannsen2},
implies  the existence of an increasing chain of equivalence relations 
$(\sim_\nu)_{1\le \nu\le
\infty}$ such that
\begin{itemize}
\item $\sim_1=\hom$;
\item if $\alpha,\beta$ are composable Chow correspondences such that 
$\alpha \sim_\mu 0$ and
$\beta\sim_\nu 0$, then $\beta\circ\alpha\sim_{\mu+\nu} 0$;
\item for any smooth projective variety $X$, there exists $\nu=\nu(X)$ 
such that
$A_{\sim_\nu}(X\times X) = A_\rat(X\times X)$.
\end{itemize}

There properties, together with the sign conjecture, imply Conjecture \ref{c3.1} by Proposition \ref{p3.2} b).
\end{proof}

\begin{rk} In fact, one has more precise but slightly weaker implications:
(Bloch-Beilinson--Murre conjecture + ``$\hom=\num$'' conjecture) $\Rightarrow$ (Voevodsky's
conjecture) $\Rightarrow$ (Kimura-O'Sul\-liv\-an conjecture [any Chow motive is
finite-dimensional]) $\Rightarrow$  (Conjecture \ref{c3.1}): see the synoptic table in
\cite[end of Ch. 12]{andre2}. 

For the first implication, see \cite[Th. 11.5.3.1]{andre2}. For the second one, see \cite[Th.
12.1.6.6]{andre2}. The third one is in Proposition \ref{p3.2} a). 
\end{rk}

\begin{defn}\label{d3.1} Let $M\in \Mot_\sim(F,A)$. For $n\in\Z$, we write $\nu(M)\ge n$ if
$M\otimes\bL^{\otimes -n}$ is effective.\footnote{By convention, we say here that a motive $N\in \Mot_\sim(F,A)$ is \emph{effective} if it is isomorphic to a motive of $\Mot_\sim^\eff(F,A)$.}
\end{defn}

\begin{prop}\label{p3.6} Suppose $A\supseteq \Q$ and the nilpotence conjecture holds for
$\sim\ge \sim'$. Then:\\ 
a) The functor $\Mot_\sim(F,A)\to \Mot_{\sim'}(F,A)$ is conservative,
and for $M\in \Mot_\sim(F,A)$, any set of orthogonal idempotents in the endomorphism ring of
$M_{\sim'}$ lifts.\\ 
b) If $M\in \Mot_\sim(F,A)$ and $M_{\sim'}$ is effective, then $M$ is
effective.\\ 
c) If $M\in \Mot_\sim(F,A)$ and $\nu(M_{\sim'})\ge n$, then $\nu(M)\ge n$.\\ d)
\cite[13.2.1]{andre2} The map $K_0(\Mot_\sim(F,A))\to K_0(\Mot_{\sim'}(F,A))$ is an isomorphism
(here, the $K_0$-groups are those of additive categories).
\end{prop}

\begin{proof} a) is classical (see \cite[Lemma 5.4]{jannsen2} for the second statement). b) By
definition, $M_{\sim'}$ effective means that $M_{\sim'}$ is isomorphic to a direct summand of
$h_{\sim'}(X)$ for some smooth projective $X$. By a), one may lift the corresponding idempotent
$e_{\sim'}$ to an idempotent endomorphism $e$ of $h_\sim(X)$, and the isomorphism
$M_{\sim'}\simeq (h_{\sim'}(X),e_{\sim'})$ to an isomorphism $M\simeq (h_\sim(X),e)$. c)
follows from b) applied to $M\otimes \bL^{\otimes -n}$. d) follows from a), since then the
functor $\Mot_\sim(F,A)\to \Mot_{\sim'}(F,A)$ is conservative and essentially surjective.
\end{proof}

The importance of Conjecture \ref{c3.1} will appear again in the next subsection and in
Section \ref{s3} (see Remark \ref{r4.3} 2) and Proposition \ref{p3.1}).

\subsection{The Chow-K\"unneth decomposition}\label{s.CK} Here we take $(A,\sim)=(\Q,\rat)$. Recall that Murre
\cite{murre2} strengthened the standard conjecture C (algebraicity of the K\"unneth projectors)
to the existence of a \emph{Chow-K\"unneth decomposition}
\[h(X)\simeq \bigoplus_{i=0}^{2d} h_i(X)\]
in $\Chow(F,\Q)$. (This is part of the Bloch-Beilinson--Murre conjecture appearing in
Proposition \ref{p3.3} b)). By Proposition \ref{p3.6} a), the nilpotence conjecture together
with the standard conjecture C imply the existence of Chow-K\"unneth decompositions.

Here are some cases where the existence of a Chow-K\"unneth decoomposition is known
independently of any conjecture:

\begin{enumerate}
\item Varieties of dimension $\le 2$ (Murre, \cite{murre}, see also \cite{scholl}). In fact,
Murre constructs for any $X$ a partial decomposition
\[h(X)\simeq h_0(X)\oplus h_1(X)\oplus h_{[2,2d-2]}(X)\oplus h_{2d-1}(X)\oplus h_{2d}(X).\]
\item Abelian varieties (Shermenev, \cite{sher}).
\item Complete intersections in $\P^N$ (see next subsection).
\item If $X$ and $Y$ have a Chow-K\"unneth decomposition, then so does $X\times Y$.
\end{enumerate}

Suppose that the nilpotence conjecture holds for $h(X)\in \Chow(F,\Q)$ and that
homological and numerical equivalences coincide on $X\times X$. The latter then implies the
standard conjecture C for $X$ \cite{kleiman}, hence the existence of a Chow-K\"unneth decomposition by the
remark above. In \cite[Th. 14.7.3 (iii)]{kmp}, it is proven:

\begin{prop}\label{p3.7}  Under these hypotheses, there exists a further decomposition  for each
$i\in [0,2d]$:
\[h_i(X)\simeq \bigoplus  h_{i,j}(X)(j)\]
such that $h_{i,j}(X)=0$ for $j\notin [0,[i/2]]$ and, for each $j$, $\nu(h_{i,j}^\hom(X))=0$
(see Definition
\ref{d3.1}). Moreover, one has isomorphisms
\begin{equation}\label{e8.4}
h_{2d-i,d-i+j}(X)\iso h_{i,j}(X)
\end{equation} 
for $i\le d$. In particular, $\nu(h_i(X))>0$ for $i>d$.
\end{prop}

Let us justify the last assertion: the isomorphisms \eqref{e8.4} imply that, when $i>d$,
$h_{i,j}(X)=0$ for $j<i-d$.

Since $\Chow^\eff(F,\Q)\to \Chow(F,\Q)$ is fully faithful, all the above (refined)
Chow-K\"unneth decompositions hold for the effective Chow motives $h(X)\in
\Chow^\eff(F,\Q)$. We deduce:

\begin{cor}\label{c3.2} Under the nilpotence conjecture and the conjecture that homological and
numerical equivalences coincide, for any smooth projective variety $X$ the image of its
Chow-K\"unneth decomposition in $\Chow^\o(F,\Q)$ is of the form
\[ h^\o(X)\simeq \bigoplus_{i=0}^d  h^\o_i(X).\]
Moreover, with the notation of Proposition
\ref{p3.7}, one has 
\[ h^\o_i(X)\simeq  h^\o_{i,0}(X) \qquad \text{for $i\le d$.}\]
\end{cor}

Examples where this conclusion is true unconditionally follow faithfully the examples where the
Chow-K\"unneth decomposition is unconditionally known:

\begin{prop} The conclusion of Corollary \ref{c3.2} holds in the following cases:
\begin{enumerate}
\item Varieties of dimension $\le 2$.
\item Abelian varieties.
\item Complete intersections in $\P^N$.
\item If $X$ and $Y$ have a Chow-K\"unneth decomposition and verify this conclusion, then so
does $X\times Y$.
\end{enumerate}
\end{prop}

\begin{proof} In cases (1) and (2), the conclusion holds because one has ``Lefschetz
isomorphisms" $h_{2d-i}(X)\iso h_i(X)(d-i)$ for $i>d$. For curves, it is trivial, for surfaces
they are constructed in
\cite{murre} (see \cite[Th. 4.4. (ii)]{scholl}: the isomorphism is constructed for $i=0,1$ and
any $X$), and for abelian varieties they are constructed in
\cite{sher}. For (3), see next subsection. Finally, (4) is clear.
\end{proof}

In the case of a surface, \cite{kmp} constructs a refined Chow-K\"unneth decomposition
\[h(X)=h_0(X)\oplus h_1(X)\oplus NS_X(1)\oplus t_2(X)\oplus h_3(X)\oplus h_4(X)\]
where $NS_X$ is the Artin motive corresponding to the Galois representation defined by
$NS(\bar X)\otimes\Q$, and $t_2(X)$ is the \emph{transcendental part of $h(X)$}. (In the
notation of Proposition \ref{p3.7}, $h_{2,0}(X)=t_2(X)$ and $h_{2,1}(X)=NS_X$.) This translates
on the birational motive of
$X$ as
\[ h^\o(X)= h^\o_0(X)\oplus  h^\o_1(X)\oplus t_2^\o(X).\]

\subsection{Motives of complete intersections} These computations will be used in 
Section \ref{s3}. Here we take $A\supseteq\Q$.

For convenience, we take the notation of \cite{deligne}: so let $X\subset \P^r$ be a smooth
complete intersection of multidegree $\underline{a}=(a_1,\dots,a_d)$, and let $n=r-d=\dim X$.
Then the cohomology of
$X$ coincides with the cohomology of $\P^{r}$ except in middle dimension \cite{deligne}, and
in particular it is fully algebraic except in middle dimension. This allows us to
easily write down a Chow-K\"unneth decomposition for $h(X)$ in the sense of
Murre \cite{murre2} (see also \cite[Cor. 5.3]{ev}):

\begin{enumerate}
\item (Murre) For each $i\ne n/2$, let $c^i\in \sZ^i(X)$ be an algebraic cycle
whose cohomology class generates $H^{2i}(X)$ (here $H$ is some Weil cohomology). Then the Chow-K\"unneth projector
$\pi_{2i}$ is given by $c^i\times c^{n-i}$. We take $\pi_j=0$ for
$j$ odd $\ne n$, and $\pi_n:= \Delta_X -\sum_{j\ne n} \pi_j$.
\item Consider the inclusion $i:X\inj \P^r$. This yields morphisms of motives
\[h(\P^r)(-d)\by{i^*} h(X)\by{i_*} h(\P^r).\]
Given the decomposition $h(\P^r)\simeq \bigoplus_{j=0}^{r} \bL^j$, this
yields for each $j\in [0,n]$ morphisms
\[\bL^j\by{i^*_j} h(X)\by{i_*^j} \bL^j\]
with composition $a=\prod a_i$. Then $(1/a) i^*_ji_*^j$ defines the
$2i$-th Chow-K\"unneth projector of $X$ (denoted $\pi_{2i}$ in (1)), except if $2i=n$. Let
$\pi_n^{prim}:=1_{h(X)} - \sum_{i=0}^{n}    (1/a) i^*_ji_*^j$: the image
$p_n(X)$ of the projector $\pi_n^{prim}$ is \emph{the primitive part} of
$h_n(X)$.
\end{enumerate}

Note that the Chow-K\"unneth projectors of (1) and (2) are actually equal. Let us record here
the corresponding (refined) Chow-K\"unneth decomposition:
\begin{equation}\label{eq3.1}
h(X)\simeq \un\oplus\bL\oplus\dots \oplus \bL^n\oplus p_n(X).
\end{equation}

\begin{lemma}\label{l3.1} a) Homological and numerical equivalences agree on all (rational)
Chow groups of $X$ provided $n$ is odd or (if $\car F=0$) the Hodge realisation of $p_n(X)$
does not contain any direct summand isomorphic to $\bL^{n/2}$.\\ b) Suppose a) is satisfied.
Then for any adequate pair $(\sim, A)$ with $A\supseteq \Q$ and any
$j\in [0,n]$, we have
\[\Mot_\sim(F,A)(\bL^j,p_n(X))=\Ker(A_j^\sim(X,A)\to A_j^\num(X,A)).\]
\end{lemma}

\begin{proof} We have
\begin{align*}
A_j^\sim(X,A)&=\Mot_\sim(F,A)(\bL^j,h(X))\\
&=\bigoplus_{i=0}^n
\Mot_\sim(F,A)(\bL^j,\bL^i) \oplus \Mot_\sim(F,A)(\bL^j,p_n(X))\\
&=\Mot_\sim(F,A)(\bL^j,\bL^j) \oplus \Mot_\sim(F,A)(\bL^j,p_n(X)).
\end{align*}

For $\sim=\hom$, we have $\Mot_\sim(F,A)(\bL^j,p_n(X))=0$ by weight reasons for
$2j\ne n$ and under the hypothesis of a) for $2j=n$ (note that the Hodge realization of
$p_n(X)$ is semi-simple, as a polarisable Hodge structure). Hence the same is true for any
$\sim$ finer than $\hom$, in particular    
$\sim=\num$. This proves a). Moreover,
$\Mot_\sim(F,A)(\bL^j,\bL^j)=A$ for any choice of $\sim$. Hence b).
\end{proof}

\eqref{eq3.1} shows that the birational motive of $X$ reduces to $\un\oplus p_n^\sim(X)^\o$. In fact,
it is possible to be much more precise:

\begin{prop}\label{p3.5} Let $\underline{a}=(a_1,\dots,a_d)$ be the multidegree of $X\subset \P^r$.\\
 a) If $a_1+\dots+a_d\le r$, $ h^\o_\rat(X) = \un$.\\
b) If $a_1+\dots+a_d > r$, $ h^\o_\num(X) \ne \un$ (equivalently, $p_n^\num(X)^\o\ne 0$)
provided $\car F=0$ or $X$ is generic.
\end{prop}

\begin{proof} a) Under the hypothesis, we conclude from Ro\v\i tman's theorem \cite{roitman}
that $CH_0(X_K) \otimes \Q =\Q$ for any extension $K/F$.\footnote{Of course we could also invoke Proposition \ref{trcc} 2) since $X$ is Fano, hence rationally chain connected, but this theorem of Campana \cite{campana} and Koll\'ar-Miyaoka-Mori \cite{kmm} was proven much later than Ro\v\i tman's work.} Assertion a) then 
follows from Proposition \ref{p3.4}. For b),
it suffices to prove the statement for homological equivalence, since the kernel of
$\Mot_\hom(F,\Q)(h(X),h(X))\to \Mot_\num(F,\Q)(h(X),h(X))$ is a nilpotent ideal (see
Propositions \ref{p3.2} b) and \ref{p3.6} a)). 

If $\car F =0$, we may use Hodge cohomology and Deligne's theorem \cite[Th. 2.5 (ii) p.
54]{deligne}. Namely, with the notation of loc. cit., the condition $p_n^\hom(X)^\o=0$ implies
$h_0^{0,n}(\underline{a})=0$, which is equivalent by loc. cit., Th. 2.5 (ii) to
\[0\le \left[\frac{n+d-\sum a_i}{\sup(a_i)} \right]\]
that is, $\sum a_i\le n+d=r$.

If $\car F>0$ and $X$ is generic, we may use Katz's theorem \cite[p. 382, Th.4.1]{katz}.
\end{proof}

\begin{rks} 1) Katz also has a result concerning a generic hyperplane section of a given complete intersection, \cite[Th. 4.2]{katz}.

2) It seems possible to remove the genericity assumption in positive characteristic 
by lifting the coefficients of the equations defining $X$ to characteristic $0$. We have not worked out the details.
\end{rks}

\section{On adjoints and idempotents}\label{s3}

We now want to examine two related questions:

\begin{enumerate}
\item Does the projection functor
$\Mot_\sim^\eff(F,A)\to
\Mot_\sim^\eff(F,A)/\sL_\sim$ have a right adjoint? This question was
raised by Luca Barbieri-Viale and is closely related to a
conjecture of Voevodsky \cite[Conj. 0.0.11]{voebei}. 
\item Is the category $\Mot_\sim^\eff(F,A)/\sL_\sim$ pseudo-abelian,  i.e., is
it superfluous to take the pseudo-abelian  envelope  in Definition \ref{d1.2}?
\end{enumerate}

The answer to both questions is ``yes" for $\sim=\num$ and $A\supseteq \Q$, as
an easy consequence of Jannsen's semi-simplicity theorem for numerical motives
\cite{jannsen}. In fact:

\begin{prop}[\protect{\cite[Prop. 7.7]{FqX}}]\label{p0} a) The projection
functor
\[\pi:\Mot_\num^\eff\to\Mot_\num^\o\] is essentially surjective.\\
b) $\pi$ has a section $i$ which is also a left and right adjoint.\\
c) The category $\Mot_\num^\eff$ is the coproduct of $\Mot_\num^\eff\otimes \bL$
and
$i(\Mot_\num^\o)$, i.e. any object of $\Mot_\num^\eff$ can be uniquely
written as a direct sum of objects of these two subcategories.
\end{prop}
In the sequel, we want to examine these questions for a general adequate pair:
see Theorems \ref{t3.1} and \ref{t3.2} for (1) and Proposition \ref{p3.1} for (2). This requires
some preparation.

\subsection{A lemma on base change}

Let $P:\sA\to \sB$ be a functor. Recall that one says that
``its" right adjoint is \emph{defined at $B\in \sB$} if the functor
\[\sA\ni A\mapsto \sB(PA,B)\]
is representable. We write $P^\sharp B$ for a representing object (unique up to
unique isomorphism).

Let
\[\begin{CD}
\sA@>\phi>> \sB\\
@V{P}VV @V{Q}VV\\
\sC@>\psi>> \sD
\end{CD}\]
be a naturally commutative diagram of pseudo-abelian additive categories, and
let $A\in \sA$.

Suppose that ``the" right adjoint $P^\sharp$ of $P$ is defined at $PA\in \sC$ and
that the right adjoint $Q^\sharp$ of $Q$ defined at
$\psi PA \simeq Q\phi A$. We then have two corresponding unit maps (adjoint to
the identities of $PA$ and $Q\phi A$)
\begin{align*}
\epsilon_P:A&\to P^\sharp P A\\
\epsilon_Q:\phi A&\to Q^\sharp Q \phi A.
\end{align*}

\begin{lemma}\label{l1.2} Suppose that $\epsilon_Q$ is an isomorphism. Then
$\phi\epsilon_P$ has a retraction. If moreover $\phi$ is full and $\Ker(\End_\sA(A)\to\End_\sB(\phi A))$ is a nilideal, then
$\epsilon_P$ has a retraction.
\end{lemma}

\begin{proof} Let $\eta_P:PP^\sharp P A\to PA$ be the counit map of the
adjunction at $PA$ (adjoint to the identity of $P^\sharp PA$), and let $u:Q\phi
 A\iso \psi P A$, $v:Q\phi
P^\sharp P A\iso \psi P P^\sharp P A$ be the natural isomorphism from $Q\phi$ to
$\psi P$ evaluated respectively at $A$ and $P^\sharp P A$. We then have a
composition
\[\begin{CD}
Q\phi P^\sharp P A@>{v}>> \psi PP^\sharp P A@>{\psi \eta_P}>>\psi PA
\end{CD}\]
which yields by adjunction a ``base change
morphism"
\[\phi P^\sharp P A\by{b} Q^\sharp\psi PA.
\]

Inspection shows that the diagram
\[\begin{CD}
\phi A@>\phi\epsilon_P>> \phi P^\sharp PA\\
@V{\epsilon_Q}VV @V{b}VV\\
Q^\sharp Q \phi A @>Q^\sharp u>> Q^\sharp \psi P A
\end{CD}\]
commutes. The first claim follows, and the second claim follows from the first.
\end{proof}

\subsection{Right adjoints}

We come back to Question (1) posed at the beginning of this section. In \cite[14.8.7]{kmp} and
\cite[7.8 3)]{FqX}, it was announced that one can show the non-existence of the right adjoint
for $\sim=\rat$, using the results of \cite[Appendix]{huber-ayoub}. The proof turns out not
to be exactly along these lines, but is closely related: see Lemma
\ref{l3.2}, Theorem \ref{t3.1} and Theorem \ref{t3.2}.

Let us abbreviate the notation to $\Mot^\eff = \Mot_\sim^\eff(F,A)$,
$\Mot^\o=\Mot_\sim^\o(F,A)$. Let $P:\Mot^\eff\to \Mot^\o$ denote the projection
functor, and let $P^\sharp$ denote its (a priori partially defined) right
adjoint. Let $\sL^\perp$ be the full subcategory of
$\Mot^\eff$ consisting of those $M$ such that $\Hom(N(1),M)=0$ for
all $N\in \Mot^\eff$. Recall from
\cite[Prop. 7.8.1]{kmp} that
\begin{itemize}
\item If $P^\sharp$ is defined at $M$, then $P^\sharp M\in \sL^\perp$;
\item The full subcategory $\Mot^\sharp$ of $\Mot^\o$ where $P^\sharp$ is
defined equals $P(\sL^\perp)$;
\item $P^\sharp$ and the restriction of $P$ to $\sL^\perp$ define
quasi-inverse equivalences of categories between $\sL^\perp$ and
$\Mot^\sharp$.
\end{itemize}

The right adjoint $P^\sharp$ is defined
at birational motives of varieties of dimension $\le 2$ for any adequate pair
$(A,\sim)$ such that $A\supseteq \Q$ by
\cite[Cor. 7.8.6]{kmp}. (The proof there is given for $(A,\sim)=(\Q,\rat)$,
but the argument works in general.)  Recall that
\[P^\sharp h^\o(C)=\un \oplus h_1(C),\qquad P^\sharp h^\o(S) = \un \oplus h_1(S) \oplus t_2(S)\]
with the notation at the end of \S \ref{s.CK}, where $C$ is a curve and $S$ is a surface. 

The following lemma gives a sufficient
condition for the nonexistence of $P^\sharp P M$ for an effective motive  $M$.

\begin{lemma}\label{l3.2} Let $(\Q,\sim)$ be an adequate pair, and let
$M\in \Mot^\eff_\sim(F,\Q)$. Assume that
\begin{thlist}
\item $M_\num\in \Mot_\num^\eff(F,\Q)$ does not contain any direct
summand divisible by $\bL$;  
\item $\Ker(\End(M)\to\End(M_\num))$ is a nilideal;
\item There exists $r>0$ such that $\Hom(\bL^r,M)\ne 0$.
\end{thlist}
Then $P^\sharp P M$ does not exist.
\end{lemma}

\begin{proof}  Suppose that $P^\sharp$ is defined at $P
M$. Consider the unit map
\begin{equation}\label{eq1.3}
\epsilon_\sim:M\to P^\sharp P M.
\end{equation}

For $\sim=\num$, $P^\sharp_\num P_\num M_\num$ exists  by Proposition \ref{p0}.
Moreover, part c) of this proposition shows that, under Condition (i) of
the lemma, $\epsilon_\num$ is an isomorphism. By Lemma
\ref{l1.2}, the image of $\epsilon_\sim$ modulo numerical equivalence then has a
retraction, and so does $\epsilon_\sim$ itself under Condition (ii). If this is
the case,
$M\in \sL^\perp$, and in particular, $\Hom(\bL^r,M)=0$ for all $r>0$, contradiction.
\end{proof}

\subsection{Counterexamples}

To give examples where the conditions of Lemma \ref{l3.2} are satisfied,
we appeal as in \cite{huber-ayoub} to the nontriviality of the Griffiths group.

We start with an example which a priori only works for a specific adequate equivalence, because
the proof is simpler. Unlike in \cite{huber-ayoub}, we don't need the full force of Clemens'
theorem \cite[Th. 0.2]{clemens}, but merely the previous results of Griffiths
\cite{griffiths}.

\begin{defn}[``Abel-Jacobi equivalence"] Let $k=\C$. For $X$ smooth projective,
$\sZ^j_\AJ(X,\Q)$ is the image of $CH^j(X)\otimes \Q$ in Deligne-Beilinson
cohomology via the (Deligne-Beilinson) cycle class map \cite{ev2}. This defines an adequate
equivalence relation.
\end{defn}

\begin{thm}\label{t3.1} Let $F=\C$ and $\sim=\AJ$. Then\\
a) Condition (ii) of
Lemma \ref{l3.2} is satisfied for any pure motive $M$.\\ 
Let $X$ be a generic hypersurface of degree $a$ in $\P^{n+1}$.\\ 
b) Condition (i) of Lemma \ref{l3.2} is
satisfied for $M=p_n(X)$ (see \eqref{eq3.1}) provided $X$ is not a quadric, a
cubic surface or an even-dimensional intersection of two quadrics, and $a\ge n+1$.\\  
c) If $n=2m-1$ is odd and
$a\ge 2+ 3/(m-1)$, then Condition (iii)  of Lemma \ref{l3.2} is satisfied for
$r=m-1$.\\
d)  $P^\sharp$ is not defined at $h^\o(X)$ in the following cases: $n$ is odd and
\begin{thlist}
\item $n=3$: $a\ge 5$.
\item $n> 3$: $a\ge n+1$.
\end{thlist}
\end{thm}

\begin{proof} a) holds because
$\Ker(\End_\AJ(M)\to \End_\hom(M))$ has square $0$ \cite[Prop. 7.10]{ev2}\footnote{A more functorial justification is: 1) Deligne-Beilinson cohomology can be computed as absolute Hodge cohomology as in \cite{absh}, 2) the category of polarisable $\Q$-mixed Hodge structures has $\Ext$-dimension $1$.} and $\Ker(\End_\hom(M)\to
\End_\num(M))$ is nilpotent.

b) By  \cite[Ex. 5 and Cor. 18]{ps}, the Hodge realisation $P_n(X)$ of $p_n(X)$
is an absolutely simple pure Hodge structure: this, together with Proposition \ref{p3.5} b), is amply sufficient to imply Condition (i) of Lemma
\ref{l3.2}.

c) By \cite[Cor. 13.2 and 14.2]{griffiths}, $\Ker(A_{m-1}^\sim(X,\Q)\to
A_{m-1}^\num(X,\Q))\ne 0$.  But by Lemma
\ref{l3.1}, this group is $\Hom(\bL^{m-1},p_n(X))$.

d) Note that, by the refined Chow-K\"unneth decomposition \eqref{eq3.1}, $P^\sharp$
is defined at $P h(X)$ if and only if it is defined at $P p_n(X)$. The conclusion now follows
from Lemma \ref{l3.2} and from collecting the results of a), b) and c).
\end{proof}

To get a conterexample with rational equivalence, we appeal to a result of Nori \cite{nori}. We
thank Srinivas for pointing out this reference.

\begin{thm}\label{t3.2} Let $X$ be a generic abelian threefold over $k=\C$. If $\sim\ge
\alg$, then $P^\sharp$ is not defined at $h^\o_\sim(X)$.
\end{thm}

\begin{proof} It is similar to that of Theorem \ref{t3.1}, except that the motive of an abelian variety is more complicated than that of a hypersurface. We only sketch the argument (details will appear elsewhere):  

It is enough to show that $P^\sharp$ is not defined at $h^\o_{3,0}(X)$, where
$h_{3,0}(X)$ is as in Proposition \ref{p3.7} (here we use that the nilpotence conjecture is true for motives of abelian varieties, see Proposition \ref{p3.2} a)). We check the conditions of Lemma \ref{l3.2} for $M=h_{3,0}(X)$. (i)
is true by definition. (ii) is true by Proposition \ref{p3.2} a). For (iii), one can show that computing
the decomposition
\[A_1^\sim(X)=\Mot_\sim^\eff(\bL,h(X))\simeq
\bigoplus_{i=0}^6\bigoplus_{j=0}^{[i/2]}\Mot_\sim^\eff(\bL,h_{i,j}(X)(j))\] 
yields a surjection
\[\Mot_\sim^\eff(\bL,h_{3,0}(X))\surj \Griff_1(X)\]
for $\sim\ge \alg$, where $\Griff_1(X)=\Ker(A_1^\alg(X)\to A_1^\num(X))$ is the Griffiths group of $X$. By Nori's theorem \cite{nori},
$\Griff_1(X)\ne 0$, and the proof is complete.
\end{proof}

\begin{rk}\label{r4.3} It is easy to get examples of any dimension $\ge 4$
by multiplying the example of Theorem \ref{t3.2} with $\P^n$.
\end{rk}

\subsection{Idempotents} We now address Question (2) from the beginning
of this section.

\begin{prop}\label{p3.1} Let $(A,\sim)$ be an adequate pair with $A\supseteq
\Q$, and let $\sM$ be a full subcategory of $\Mot_\sim^\eff(F,A)$ closed under
direct summands. If the nilpotence conjecture
\ref{c3.1} holds for the objects of $\sM$, then the category
$\sM/\sL_\sim$ is pseudo-abel\-ian.
\end{prop}

\begin{proof} Let $\sM_\num$ denote the pseudo-abelian envelope of the image of
$\sM$ in $\Mot_\num^\eff(F,A)$. We have a commutative diagram of categories:
\[\begin{CD}
\sM@>P>> \sM/\sL_\sim\\
@V{\pi}VV @V{\bar\pi}VV \\
\sM_\num@>P_\num>> \sM_\num/\sL_\num
\end{CD}\]

Under the hypothesis,
$\pi$ is essentially
surjective (one can lift idempotents).
Hence $\bar\pi$  is essentially surjective as well.  Since $P$ is essentially
surjective and $\pi,P_\num$ are full,
$\bar\pi$ is full, and its kernel is locally nilpotent as a quotient of the
kernel of
$\pi$ (fullness of $P$). Thus 
$\bar\pi$ is full, essentially surjective and conservative. 

Since $\Mot_\num^\eff(F,A)$ is abelian semi-simple,  $\sM_\num$ is also abelian 
semi-simple, hence so is $\sM_\num/\sL_\num$ which is in particular
pseudo-abelian. 

Let now
$M\in
\sM/\sL_\sim$, and let
$p=p^2\in \End(M)$. Write $M_\num\simeq M_1\oplus M_2$, where $M_1=\IM p_\num$
and $M_2=\Ker p_\num$. By essential surjectivity, we may lift $M_1$ and $M_2$ to
objects $\tilde M_1,\tilde M_2\in \sM/\sL_\sim$.

By fullness, we may lift the isomorphism $M_1\oplus M_2\iso M_\num$ to a
morphism $\tilde M_1\oplus \tilde M_2\to M$ in
$\sM/\sL_\sim$, and this lift is an
isomorphism by conservativity. This concludes the proof.
\end{proof}

\begin{ex}\label{e3.1} Proposition \ref{p3.1} applies taking for $\sM$ the
category of motives of abelian type (direct summands of the tensor product of an
Artin motive and the motive of an abelian variety), since such motives are
finite-dimensional (Kimura \cite{kimura}).
\end{ex}

The situation when $A$ does not contain $\Q$, for example
$A=\Z$, is unclear.

\section{Birational motives and birational categories}\label{s4}

In this section, we relate the categories studied in \cite{Birat} with the
categories of pure birational motives introduced here.

\subsection{}
 From \eqref{eq7.4}, we get a composite
functor:
\begin{equation}\label{eq5.1}
 S_r^{-1}\Sm^\proj(F) \to S_r^{-1}\Chow^\eff(F)\to \Chow^\o(F).
\end{equation}

The morphisms in the first category can be
described by means of
$R$-equivalence classes  \cite[Th. 6.6.3, Cor. 6.6.4 and Rk. 6.6.5]{Birat}; by Lemma
\ref{ehomrat}, those in the last category can be described by means of Chow
groups of $0$-cycles. One checks easily that the action of the composite functor
on Hom sets is just the map which sends $R$-equivalence classes of
rational points to $0$-cycles modulo rational equivalence. This puts this
map within a functorial setting.

Let us now recall further results from \cite{Birat}. Let  $\place(F)$ denote the category of
finitely generated extensions of $F$, with $F$-places as morphisms. In \cite[(4.3)]{Birat}, we constructed a functor
\[\place_*(F)^\op\to S_b^{-1} \Sm^\proper(F)\]
hence a functor
\[S_r^{-1}\place_*(F)^\op\to S_r^{-1} \Sm^\proper(F)\]
where $\place_*(F)$ denotes the full subcategory of $\place(F)$ defined
by those $K/F$ which have a cofinal set of smooth proper models, and $S_r\subset Ar(\place(F))$ denotes the set of purely transcendental extensions. The same arguments as in loc. cit. give an analogous functor
\begin{equation}\label{5.1sharp}S_r^{-1}\place_\sharp(F)^\op\to S_r^{-1} \Sm^\proj(F)
\end{equation}
where $\place_\sharp(F)$ has the same definition as $\place_*(F)$, replacing ``smooth proper'' by ``smooth projective''. Composing \eqref{5.1sharp} with \eqref{eq5.1}, we get a functor
\begin{equation}\label{eq5.3}
S_r^{-1}\place_\sharp(F)^\op\to\Chow^\o(F).
\end{equation}

We can describe the 
image under this functor of a place $\lambda:K\tto
L$ in $CH_0(X_L)$, where $X$ is a smooth projective model of $K$: it is just
the class of the centre of $\lambda$. Hence the image of \eqref{eq5.3}
on morphisms consists of the classes of $L$-rational points. This answers a
question of D\'eglise.

In characteristic $0$, 
$\place_\sharp(F)=\place(F)$ by resolution of singularities and $S_r^{-1}\Sm^\proj(F)\iso S_r^{-1}\Sm(F)$ by \cite[Prop. 8.5]{localisation}. In 
characteristic $p$, we would 
ideally like to get functors
\begin{align*}
S_r^{-1}\place(F)^{op}&\to \Chow^\o(F)\\
S_r^{-1}\Sm(F)&\to \Chow^\o(F)
\end{align*}
extending \eqref{eq5.1} and \eqref{eq5.3}.   Constructing the first functor looks technically 
difficult: we 
shall content ourselves with extending \cite[Rk. 7.4]{FqX} to all 
finitely generated fields $K/F$, by using an adjunction result from  \cite{adjoints}; this will not be used in the rest of the paper. The second functor is constructed in \cite[Corollary 2.4.2]{birat-tri}.

\begin{prop}\label{ph0} Let $p$ be the exponential characteristic of $F$.\\
a) There is a unique functor 
(up to unique isomorphism)  
\[h^\o:S_r^{-1}\field(F)^{op}\to \Chow^\o(F,\Z[1/p])\]
such that, for any $K\in \field(F)$ and any $Y\in\Sm^\proj(F)$, one has 
\begin{equation}\label{eq5.2}
\Chow^\o(F,\Z[1/p])(h^\o(K),h^\o(Y)) \simeq CH_0(Y_K)\otimes\Z[1/p].
\end{equation}
This functor transforms purely inseparable extensions into 
isomorphisms.\\
b) If $K\subseteq L$, the map $h^\o(L)\to h^\o(K)$ has a section.\\
c) We have $h^\o(K) =h^\o(X)$ if $K = F(X)$ for a smooth projective variety 
$X$. Moreover, if $K=F(X)$, $L=F(Y)$ with $X,Y$ smooth projective, and if 
$f:K\to L$ 
corresponds to a rational map $\phi:Y\ttto X$, then $h^\o(f)$ is given by 
the graph of $\phi$.
\end{prop}

\begin{proof} a) Note that the isomorphism \eqref{eq5.2} determines 
$h^\o(K)$ up to unique isomorphism, by Yoneda's lemma. By Lemma
\ref{ehomrat} applied over $K$, this isomorphism may 
be rewritten as
\[\Chow^\o(F,\Z[1/p])(h^\o(K),h^\o(Y)) \simeq \Chow^\o(K,\Z[1/p])(\un_K,h^\o(Y_K)).\]
where $\un_K=h^\o(\Spec K)$ is the unit object of $\Chow^\o(K,\Z[1/p])$.

By \cite[Th. 6.5]{adjoints}, the base-change functor
\[\Chow^\o(F,\Z[1/p])\to \Chow^\o(K,\Z[1/p])\]
has a left adjoint $l_{K/F}$. Therefore we may define 
$h^\o(K)=l_{K/F}(\un_K)$. 

Suppose $F\to K\by{f} L$ are successive finitely generated 
extensions. Since the base-change of $\un_K$ is $\un_L$, the identity map 
$\un_L \to \un_L$ gives by adjunction a map
\[l_{L/K}\un_L\to \un_K\]
hence a map
\[h^\o(f):h^\o(L)=r_{L/F}(\un_L)\to r_{K/F}(\un_K)= h^\o(K).\]

We just used the transitivity of adjoints; using it a second time on a 
$3$-layer extension shows that we have indeed defined a functor 
$\field(F)^{op}\to \Chow^\o(F,\Z[1/p])$.

Suppose that $L=K(t)$. Then $l_{L/K}(\un_L) = h^\o(\P^1)=\un_K$, hence 
$h^\o(f)$ is an isomorphism. This shows that our functor induces a functor 
$h^\o:S_r^{-1}\field(F)^{op}\to \Chow^\o(F,\Z[1/p])$, as required.

Suppose now that $K\by{f} L$ is a finite and purely inseparable 
extension of finitely generated fields over $F$. If 
$X$ is a smooth projective $K$-variety, the map $CH_0(X)\otimes\Z[1/p]\to 
CH_0(X_L)\otimes\Z[1/p]$ is an isomorphism by Lemma \ref{l4.1a}: this shows 
that $l_{L/K}(\un_L)=\un_K$, hence that $h^\o(f)$ is invertible.

b) The proof is the same as in \cite[Rk. 7.4]{FqX}: write $L$ as a finite 
purely inseparable  
extension of a finite separable extension of a purely transcendental 
extension of $K$. Then a) reduces us to the case where $L/K$ is finite 
and separable. We may write $L=\Spec X$ where $X$ is a 
$0$-dimensional smooth projective $K$-variety, and 
$l_{L/K}(\un_L)=h^\o(X)$. 
The conclusion now follows from Lemma \ref{l7.1}.

c) If $K=F(X)$ 
for $X$ smooth projective, then Lemma \ref{ehomrat} and Yoneda's lemma 
show that $h^\o(K) \simeq
h^\o(X)$. 
For the claim on morphisms, we are reduced (again by Yoneda's lemma) to 
determining the map
\[\Chow^\o(F,\Z[1/p])(h^\o(K),h^\o(Z))\by{h^\o(f)^*}\Chow^\o(F,\Z[1/p])(h^\o(L),h^\o(Z))
\]
for a smooth projective $F$-variety $Z$. By definition 
of $h^\o(f)$, an adjunction computation shows that this map may be 
rewritten as the map
\begin{multline*}
CH_0(Z_K)\otimes\Z[1/p]
=\Chow^\o(K,\Z[1/p])(\un_K,h^\o(Z_K))\\
\to\Chow^\o(L,\Z[1/p])(\un_L,h^\o(Z_L)) 
=CH_0(Z_L)\otimes\Z[1/p]
\end{multline*}
given by extension of scalars. The conclusion now follows from Lemma 
\ref{l2.3}.
\end{proof}

\section{Birational motives and cycle modules}\label{s.6}

In \cite{rost}, Rost introduced the notion of cycle module and cycle cohomology; he proved in \loccit, Cor. 12.10 that for any cycle module $M$, $A^0(X,M)$ is a birational invariant of smooth projective varieties $X$. In \cite[Cor. 3.5]{merk}, he extended this to $A_0(X,M)$ by introducing the category $\Chow^\o(F)$ of Definition \ref{d1.2a} (independently from this paper). In the first subsection, we essentially reproduce \S 3 of \cite{merk}; we don't claim any originality here, but hope this will be a service to the reader since this preprint remains unpublished. In the second subsection, we connect these results with more recent work of Merkurjev.

To lighten notation, we drop the reference to the base field $F$ in the relevant categories.

\subsection{The functors $A^0$ and $A_0$} Let $M=(M_n)_{n\in\Z}$ be a cycle module over $F$ in the sense of Rost \cite{rost}: recall that this is a functor from $\field$ to graded abelian groups, provided with extra structure (transfers, residues, cup-products by units) subject to certain axioms. To a smooth variety $X\in \Sm$, one associates its \emph{cycle cohomology with coefficients in $M$} \cite[\S 5]{rost}
\[A^p(X,M_n) = H(\dots\by{\partial} \bigoplus_{x\in X^{(p)}} M_{n-p}(F(x))\by{\partial} \dots)\]
where the differentials $\partial$ are induced by the residue homomorphisms. We also have the homological notation
\[A_p(X,M_n) = H(\dots\by{\partial} \bigoplus_{x\in X_{(p)}} M_{n+p}(F(x))\by{\partial} \dots)\]
so that $A_p(X,M_n)=A^{d-p}(X,M_{d+n})$ if $X$ is purely of dimension $d$.

\begin{prop}\label{p6.3} a) Let $X,Y$ be two smooth projective varieties and let $\alpha\in CH_{\dim X}(X\times Y)$ be a Chow correspondence. Then $\alpha$ induces homomorphisms
\[\alpha^*:A^p(Y,M_n)\to A^p(X,M_n), \quad \alpha_*:A_p(X,M_n)\to A_p(Y,M_n)\]
which make $A^p(-,M_n)$ (resp. $A_p(-,M_n)$) a contravariant (resp. covariant) functor on $\Chow^\eff$.\\
b) Suppose that $\alpha\in\sI_\rat(X,Y)$, where $\sI_\rat$ is as in Proposition \ref{l4.2'}. Then $\alpha^*A^0(Y,M_n)=0$ (resp. $\alpha_* A_0(X,M_n)=0$).
\end{prop}

\begin{proof} a) follows easily from the functoriality of cycle cohomology \cite[Prop. 4.6, \S 13, \S 14]{rost}. Namely, we define $\alpha^*$ as the composition
\begin{multline}\label{eq.comp}
A^p(Y,M_n)\by{p_Y^*}A^p(X\times Y,M_n)\\
\by{\cup \alpha} A^{p+\dim Y}(X\times Y,M_{n+\dim Y})\by{p_{X*}} A^p(X,M_n)
\end{multline}
where $\cup\alpha$ is cup-product with $\alpha$ as in \cite[\S 14]{rost}, and $\alpha_*$ similarly. Checking the identities $(\beta\circ \alpha)^*=\alpha^*\circ \beta^*$ and $(\beta\circ \alpha)_* = \beta_*\circ \alpha_*$ is a routine matter, using the compatibility of cup-product with pull-backs and the projection formula (ibid.).

To prove b), we may assume $X$ irreducible; let $Z\subset X$ be a proper closed subset such that $\alpha$ is supported on $Z\times Y$, and let  $U=X-Z$. We consider the cases of $\alpha^*$ and $\alpha_*$ separately. 

In the first case,  we observe that \eqref{eq.comp} also makes sense for $X$ smooth (not necessarily projective) and that $A^0(X,M_n)\to A^0(U,M_n)$ is injective (both groups being subsets of   $M_n(F(X))$). Therefore it suffices to see that \eqref{eq.comp} is $0$ when $X$ is replaced by $U$, which is obvious since $\alpha_{|CH_{\dim X}(U\times Y)}=0$.

In the second case, we generalise the argument in the proof of Proposition \ref{l4.2'}: if $x\in X_{(0)}$, it suffices to show that the composition
\begin{multline*}
M_n(F(x))\by{i_{x*}} A_0(X,M_n)=A^{\dim X}(X,M_{n+\dim X})\\
\by{p_Y^*}A^{\dim X}(X\times Y,M_{n+\dim X})
\by{\cup \alpha} A^{\dim X+\dim Y}(X\times Y,M_{n+\dim X+\dim Y})\\
\by{p_{Y*}} A^{\dim Y}(Y,M_{n+\dim Y})=A_0(Y,M_n)
\end{multline*}
is $0$. If $q_Y:x\times Y\to x$ is the first projection, we have
\[p_Y^*i_{x*} = (i_x\times 1_Y)_*q_Y^*\]
\cite[Prop. 4.1 (3)]{rost}. For $a\in M_n(F(x))$, we now have
\[
p_Y^*i_{x*} a\cup \alpha = (i_x\times 1_Y)_*q_Y^*a\cup \alpha 
= (i_x\times 1_Y)_*(q_Y^*a\cup (i_x\times 1_Y)^*\alpha)
\]
by the projection formula \cite[14.5]{rost}. As in the proof of Proposition \ref{l4.2'} we reduce to the case where $x\notin Z$, and then $(i_x\times 1_Y)^*\alpha=0$.
\end{proof}

From Proposition \ref{p6.3} b), we immediately deduce:

\begin{cor}\label{c6.1} a) For any cycle module $M$ and any $n\in \Z$, the assignment
\[\Sm^\proj\ni X\mapsto A^0(X,M_n) \text{ (resp. } A_0(X,M_n))\]
extends to a contravariant (resp. a covariant) additive functor
\[A^0(-,M_n) \text{ (resp. } A_0(-,M_n)):\Chow^\o\to \Ab.\]
b) Let $X\in \Sm^\proj$ be such that $h^\o(X)\simeq \un\in \Chow^\o(F)$. Then the maps
\[M_n(F)\to A^0(X,M_n), \quad A_0(X,M_n)\to M_n(F)\]
induced by the structural map $\pi_X:X\to \Spec F$ are isomorphisms for any cycle module $M$ and any $n\in \Z$.\qed
\end{cor}

This proves the implication (iv) $\Rightarrow$ (v) in Proposition \ref{p3.4}.

\subsection{Relationship with Merkurjev's work} For $A^0(X,M_n)$, Corollary \ref{c6.1} b) is part of a theorem of Merkurjev : 
\begin{prop}[\protect{\cite[Th. 2.11, (3) $\Rightarrow$ (1)]{merk2}}]\label{q6.1}\it
If $CH_0(X_E)\iso \Z$ for any extension $E/F$, then $M_n(F)\iso A^0(X,M_n)$ for all cycle modules $M$ and all $n\in \Z$. 
\end{prop}

Indeed, this condition is equivalent to $h^\o(X)\simeq \un$ in $\Chow^\o$ by (iv) $\iff$ (i) in Proposition \ref{p3.4}.

Merkurjev proves the converse implication. For this, he defines a  cycle module $K^X$ such that
\[K^X_n(E)=A_0(X_E,K_n)\]
for any extension $E/F$. Here, $K$ is the cycle module given by Milnor $K$-theory. He shows:

\begin{thm}[\protect{\cite[Th. 2.10]{merk2}}] \label{t6.1}The functor
\begin{align*}
\CM &\to \Ab\\
M&\mapsto A^0(X,M_0)
\end{align*}
from the category of cycle modules to abelian groups is corepresented by $K^X$.
\end{thm}

See \cite[Th. 1.3]{modcyclnr} for a generalisation to non-proper $X$'s.

 Let us give a proof of the converse to Proposition \ref{q6.1} via birational motives, using only the existence of $K^X$ and thus completing the proof of Proposition \ref{p3.4}. Let us say that a cycle module $M$ is \emph{connected} if $M_n=0$ for $n<0$: we note that
\begin{equation}\label{eq6.3}
A^0(X,M_0)=M_0(F(X))\quad \text{if $M$ is connected.}
\end{equation}
 
As  $K^X$ is connected and $K^X_0(E) = CH_0(X_E)$, the condition $K_0^X(F)\allowbreak\iso A^0(X,K_0^X)$ translates as 
$CH_0(X)\iso CH_0(X_{F(X)})$, which in turn implies Condition (iii) in Proposition \ref{p3.4}.

We are now going to use Theorem \ref{t6.1} to clarify the relationship between birational motives and cycle modules.

\begin{thm} Let $\Mod\Chow^\o$ be the category of of additive
contravariant functors from $\Chow^\o$ to $\Ab$. The functor
\[A^0:\CM\to \Mod\Chow^\o\]
from Corollary \ref{c6.1} a) has a fully faithful left adjoint $\Lambda\mapsto K^\Lambda$; the essential image of this left adjoint is contained in the full subcategory of connected cycle modules.
\end{thm}

\enlargethispage*{20pt}

\begin{proof}  We first observe that $X\mapsto K^X$ extends to a functor 
\[\Chow^\o\to \CM\]
thanks to Corollary \ref{c6.1} a) (case of $A_0$). 
Let $\Lambda\in\Mod\Chow^\o$. We define
\[K^\Lambda=\colim_{y(X)\to \Lambda} K^X\]
where $y:\Chow^\o\to \Mod\Chow^\o$ is the additive Yoneda functor, and the colimit is taken on the comma category $y\downarrow \Lambda$ \cite[Ch. II, \S 6]{maclane}. Since $K^X$ is connected for any smooth projective $X$, $K^\Lambda$ is connected. For a cycle module $M$, the identity
\[\CM(K^\Lambda,M)\simeq  \Mod\Chow^\o(\Lambda,A^0(M))\]
follows from Theorem \ref{t6.1} and Yoneda's lemma, thus proving the existence of the left adjoint and the statement on its essential image.

It remains to show that $\Lambda\mapsto K^\Lambda$ is fully faithful or, equivalently, that the unit map
\[\Lambda\to A^0(K^\Lambda)\]
is an isomorphism for all $\Lambda$. Let $Y\in \Sm^\proj$: we need to show that
\[\Lambda(h^\o(Y))\to A^0(Y,K_0^\Lambda)=K_0^\Lambda(F(Y))\]
is an isomorphism, where we just used \eqref{eq6.3}. We compute:
\begin{multline*}K_0^\Lambda(F(Y))=\colim_{y(X)\to \Lambda}K_0^X(F(Y))=\colim_{y(X)\to \Lambda}CH_0(X_{F(Y)})\\
= \colim_{y(X)\to \Lambda}\Chow^\o(h^\o(Y),h^\o(X))\\
= \colim_{y(X)\to \Lambda}y(h^\o(X))(h^\o(Y))=\Lambda(h^\o(Y)).
\end{multline*}
\end{proof} 

We come back to the essential image of the functor $K^?$ in \cite[\S 4.2]{birat-tri}.

\section{Locally abelian schemes}\label{s5}

In this section, $F$ is perfect. We drop it from the notation for relevant categories.

\subsection{The Albanese scheme of a smooth projective variety}

\begin{defn}\label{d4.1} a) Let $X$ be a smooth separated $F$-scheme (not
necessarily of finite type). For each connected
component
$X_i$ of $X$, let $E_i$ be its field of constants, that is, the algebraic
closure of $F$ in $F(X_i)$. We define
\[\pi_0(X)=\coprod_i \Spec E_i.\]
There is a canonical $F$-morphism $X\to \pi_0(X)$; $\pi_0(X)$ is called
the \emph{scheme of constants} of $X$.\\
b) If $\dim X=0$ (equivalently $X\iso \pi_0(X)$), we write $\Z[X]$ for the
$0$-dimensional group scheme
representing the \'etale sheaf $f_*\Z$, where $f:X\to \Spec F$ is the
structural morphism.
\end{defn}

\begin{defn} a) For an $F$-group scheme $G$, we denote by $G^0$ the kernel
of the canonical map $G\to \pi_0(G)$ of Definition \ref{d4.1}: this is the
\emph{neutral component} of $G$.\\
b) An $F$-group scheme $G$ is called a \emph{lattice} if $G^0=\{1\}$ and
the geometric fibre of $\pi_0(G)(=G)$ is a free finitely generated abelian
group. 
\end{defn} 

\begin{defn}[\protect{\cite{ram}}] a) Recall that a \emph{semi-abelian variety} is an extension of an abelian variety by a 
torus. We denote by $\SAb$ the category of semi-abelian $F$-varieties, and by $\Ab$ the full subcategory of abelian 
varieties.\\
b) We denote by $\SAbS$ the full subcategory of the category
of commutative $F$-group schemes consisting of those objects $\sA$ such that
\begin{itemize}
\item $\pi_0(\sA)$ is a lattice;
\item $\sA^0$ is a semi-abelian variety.
\end{itemize}
Objects of $\SAbS$ will be called \emph{locally semi-abelian
$F$-schemes}.\\
c) We denote by $\AbS$ the full subcategory of $\SAbS$ consisting of those $\sA$ such that $\sA^0$ is an abelian 
variety. Its objects are called \emph{locally abelian $F$-schemes}.
\end{defn}

Note that $\SAbS$ is a Serre subcategory of the abelian category of commutative $F$-group schemes locally of finite type (\cf \cite[Exp. VI, Prop. 5.4.1 and Th. 5.4.2]{sga3}; in particular it is abelian, and $\AbS$ is idempotent-closed in $\SAbS$, hence pseudo-abelian.

For any smooth $F$-variety $X$, let $\sA_{X/F}=\sA_X$ be the
Albanese scheme of $X$ over $F$ \cite{ram}: it is an object of $\SAbS$ and there is a canonical morphism
\begin{equation}\label{eq6.2}
\phi_X:X\to \sA_X
\end{equation}
which is universal for morphisms from $X$ to objects of $\SAbS$. There is an exact sequence of
group schemes
\[0\to \sA^0_X\to \sA_X\to\Z[\pi_0(X)]\to 0\] 
where $\sA^0_X$ is the Albanese variety of $X$ (a semi-abelian
variety) and $\pi_0(X)$ has been defined above.

The aim of this section is to endow $\SAbS$ and $\AbS$ with  
symmetric monoidal structures, and to relate the latter one to birational motives (see Propositions 
\ref{p4.3} and \ref{p6.1}).

Let us recall from \cite{ram} a description of $\sA_X$. Let $\Z[X]$ be the
``free" presheaf on $F$-schemes defined by $\Z[X](Y)=\Z[X(Y)]$ and
$\sZ_{X/F}=\sZ_X$ the
associated sheaf on the big fppf site of $\Spec F$. Then $\sA_X$ is the
universal representable quotient of $\sZ_X$. In other words, there is a
homomorphism
\[\sZ_X\to \sA_X\]
where $\sA_X$ is considered as a representable sheaf, which is universal for
homomorphisms from $\sZ_X$ to sheaves of abelian groups representable by a 
locally semi-abelian $F$-scheme.

\enlargethispage*{20pt}

Let us also denote by $P_X$ the universal torsor under $\sA_X^0$
constructed by Serre \cite{serre}. There is a map 
$X\by{\tilde\phi_X} P_X$
which is universal for maps from $X$ to torsors under semi-abelian varieties.
The torsor
$P_X$ and the group scheme $\sA_X$ have the same class in
$\Ext^1_{(Sch/F)_\et}(\pi_0(\sA_X),\sA_X^0)\allowbreak=H^1_\et(\pi_0(X),\sA_X^0)$
(here we identify $\sA_X^0$ with the
corresponding representable \'etale sheaf over the big \'etale site of
$\Spec F$). A beautiful concrete description of this correspondence is
given in \cite[1.2]{ram}. The map $\tilde\phi_X$ induces an isomorphism
\[\sA_X\iso \sA_{P_X}.\]

We repeat some properties of $\sA_X$ as taken from \cite[Prop. 1.6
and Cor. 1.12]{ram} and add one.

\begin{prop}\label{p4.1} a) $\sA_X$ is covariant in $X$.\\ 
b) Let $K/F$ be an extension. Then the natural map 
\[\sA_{X_K/K}\to\sA_{X/F}\otimes_F K\] 
stemming from the universal property is an isomorphism.\\
c) If $X=Y\coprod Z$, then the natural map $\sA_{Y/F}\oplus\sA_{Z/F}\to
\sA_{X/F}$ is an isomorphism.\\
d) Let $E/F$ be a finite extension. For any $E$-scheme $S$, let
$S_{(F)}$ denote the (ordinary) restriction of scalars of $S$, i.e. we
view $S$ as an $F$-scheme.  Then there is a
natural isomorphism for $X$ smooth
\[R_{E/F}\sA_{X/E}\iso \sA_{X_{(F)}/F}\]
where $R_{E/F}$ denotes Weil's restriction of scalars.
\end{prop}

\begin{proof} The only thing which is not in \cite{ram} is d). We shall construct
the isomorphism by descent from c), using b).

Let $f:\Spec E\to \Spec F$ be the structural morphism. Recall that, for any
abelian sheaf $\sG$ on $(Sch/E)_\et$, the trace map defines an isomorphism
\cite[Ch. V, Lemma 1.12]{milne2}
\[f_*\sG\iso f_!\sG\]
where $f_!$ (\resp $f_*$) is the left (\resp right) adjoint of the
restriction functor $f^*$. This isomorphism is natural in $\sG$.

This being said, the additive version of Yoneda's lemma immediately yields
\[f_!\sZ_{X/E}=\sZ_{X_{(F)}/F}\]
hence a composition of homomorphisms of sheaves
\begin{equation}\label{eq4.1}f_*\sZ_{X/E}\iso\sZ_{X_{(F)}/F}\to
Shv(\sA_{X_{(F)}/F})
\end{equation}
where, for clarity, $Shv(\sA_{X_{(F)}/F})$ denotes the sheaf associated to
the group scheme $\sA_{X_{(F)}/F}$. We also have a chain of homomorphisms
\begin{equation}\label{eq4.2}
f_*\sZ_{X/E}\to f_*Shv(\sA_{X/E})\iso Shv(R_{E/F}\sA_{X/E})
\end{equation}
where the last isomorphism is formal. If we can prove that \eqref{eq4.1}
factors through \eqref{eq4.2} into an isomorphism, we are done by Yoneda.

In order to do this, we may assume via b) that $F$ is algebraically closed,
hence that $f$ is completely split. Then the claim follows from c).\end{proof}

We record here similar properties for the torsor $P_X=P_{X/F}$ (proofs are
similar):

\enlargethispage*{20pt}

\begin{prop} a) $X\mapsto P_X$ is a functor.\\
b) Let $K/F$ be an extension. Then the natural map $P_{X_K/K}\to
P_{X/F}\otimes_F K$ stemming from the universal property is an
isomorphism.\\
c) If $X=Y\coprod Z$, then there is an isomorphism 
$P_{Y/F}\times P_{Z/F}\iso
P_{X/F}$ which is natural in $(Y,Z)$.\\
d) Let $E/F$ be a finite extension.  Then there is a
natural isomorphism
\[P_{X_{(F)}/F}\to R_{E/F} P_{X/E}.\Box\]
\end{prop}

(In c), the map stems from the fact that coproducts correspond to
scheme-theoretic products in an appropriate category of torsors.)

\subsection{The tensor category of locally semi-abelian schemes}

Recall the Yoneda full embedding $Shv:\SAbS\to Ab((Sch/F)_\et)$,
where the latter is the category of sheaves of abelian groups over the big
\'etale site of $\Spec F$. 

\begin{lemma} \label{l4.1} a) If a sheaf $\sF\in  Ab((Sch/F)_\et)$ is an 
extension of a lattice $L$ by a semi-abelian variety $A$, 
it is represented by an object of $\SAbS$. \\
b) Let $A$ be a semi-abelian variety and $L$ a
lattice. Then the
\'etale sheaf $B=A\otimes L$ is represented by a semi-abelian variety.
\end{lemma}

\begin{proof} a) If $L$ is constant, then the choice of a basis of $L$ determines
a section of the projection $\sF\to Shv(L)$, hence an isomorphism
$\sF\simeq Shv(A)\oplus Shv(L)$. Then $\sF$ is represented by
$\coprod_{l\in L} A$. In general, $L$ becomes
constant on some finite extension $E/F$, hence $\sF_E$ is representable.
By full faithfulness, the descent data of $\sF_E$ are morphisms of
schemes; then we may apply \cite[Cor. V.4.2 a) or b)]{serre2}.

b) Same method as in a).\end{proof}

\begin{ex} If $L=\Z[\Spec E]$, where $E$ is an \'etale $F$-algebra, then
$A\otimes L=R_{E/F} A_E$.
\end{ex}

We shall also need:

\begin{lemma}\label{cB.1} Let $F$ be a field, $G_1,G_2,G_3$ be three semi-abelian $F$-varieties, and let $\phi:G_1\times 
G_2\to G_3$
be an $F$-morphism. Assume that $\phi(g_1,0)=\phi(0,g_2)=0$ identically. Then $\phi=0$.
\end{lemma}

\begin{proof} By \cite[Lemma 3]{torenis}, $\phi$ is a homomorphism and the conclusion is obvious.
\end{proof}

Let $\sA,\sB\in \SAbS$. Viewing them as \'etale sheaves, we may consider
their tensor product $\sA\otimes_{shv}\sB$. This tensor product contains
the subsheaf $\sA^0\otimes_{shv}\sB^0$, which is clearly not representable.
We
define
\[\sA\otimes_\rep\sB=\sA\otimes_{shv}\sB/\sA^0\otimes_{shv}\sB^0.\]

\begin{prop}\label{p4.2} a) $\sA\otimes_\rep\sB$ is representable by an
object of 
$\SAbS$.\\
b) For $X,Y\in \Sm$, the natural map
\[\sZ_X\otimes_{shv}\sZ_Y=\sZ_{X\times Y}\to \sA_{X\times Y}\]
factors into an
isomorphism 
\[\sA_X\otimes_\rep\sA_Y\iso \sA_{X\times Y}.\]
\end{prop}

(This corrects \cite[Cor. 1.12 (vi)]{ram}.)

\begin{proof}  a) We have a short exact sequence
\[0\to \sA^0\otimes \pi_0(\sB)\oplus \sB^0\otimes\pi_0(\sA)\to
\sA\otimes_\rep\sB\to \pi_0(\sA)\otimes\pi_0(\sB)\to 0.\]

By Lemma \ref{l4.1} b), the left hand side is representable by a
semi-abelian variety, and the right hand side is clearly a lattice. We
conclude by Lemma \ref{l4.1} a).

b) It is enough to show that this holds over the algebraic closure of $F$.
Using Proposition \ref{p4.1} c) (and the similar statement for
$\sZ$), we may assume that $X$ and $Y$ are connected. We
shall show more generally that, for any locally semi-abelian scheme $\sB$
and any map
$X\times Y\to \sB$, the induced sheaf-theoretic map
\begin{equation}\label{eq6.1}
\sZ_X\otimes_{shv}\sZ_Y\to \sB\end{equation} 
factors through
$\sA_X\otimes_\rep\sA_Y$. By a), this will show that the latter has the
universal property of
$\sA_{X\times Y}$.

For $n\in\Z$, we denote by $\sZ_X^n$ or $\sA_X^n$ the inverse
image of $n$ under the augmentation map $\sZ_X\to \Z$
or $\sA_X\to \Z$ stemming from the structural morphism $X\to \Spec F$. It
is a subsheaf of $\sZ_X$ or $\sA_X$, and $\sA_X^n$ is clearly
representable (by a variety $\bar F$-isomorphic to the semi-abelian variety
$\sA_X^0$). We shall also identify varieties with representable sheaves:
this should create no confusion in view of Yoneda's lemma.

We first show that \eqref{eq6.1} factors through
$\sA_X\otimes_{shv}\sA_Y$. It suffices to show that the composition
\[\sZ_X\times Y\to \sZ_X\otimes_{shv}\sZ_Y\to \sB\]
factors through $\sA_X\times Y$, and to conclude by symmetry. But
$X\times Y$ is connected, so its image in $\sB$ falls in some connected
component
$\sB^t$
of $\sB$, which is a torsor under $\sB^0$; applying the ``Variation en
fonction d'un param\`etre" statement in \cite[p. 10-05]{serre}, we see
that it extends to a morphism
$\sA_X^1\times Y\to
\sB^t$. Including
$\sB^t$ into $\sB$, we get a
commutative diagram
\[\begin{CD}
\sA_X^1\times Y@>>>\sB\\
@A{}AA @A{}AA\\
\sZ_X^1\times Y@>>> \sZ_X\times Y.
\end{CD}\]

Let $\sK=\Ker(\sZ_X\to \sA_X)=\Ker(\sZ_X^0\to \sA_X^0)$. The above diagram shows
that the following diagram
\[\begin{CD}
\sK\times\sZ^1_X\times Y@>a>> \sZ_X^1\times Y\\
@V{c}VV @V{d}VV\\
\sZ^1_X\times Y@>b>> \sB
\end{CD}\]
commutes, where $a$ is given by the action of $\sK$ on
$\sZ^1_X$ by left translation and $c$ is given by
$(k,z,y)\mapsto (z,y)$. Since $b$ is a homomorphism in
the first variable, this implies the desired factorisation.

We now show that the composition 
\[\sA^0_X\otimes_{shv}\sA^0_Y\to
\sA_X\otimes_{shv}\sA_Y\to \sB\] 
is $0$. It is sufficient to show that the
composition of this map with the inclusion $\sA^0_X\times\sA^0_Y\to
\sA^0_X\otimes\sA^0_Y$ is $0$. But $\sA^0_X\times\sA^0_Y$ is connected,
hence its image falls in some connected component, in fact in $\sB^0$.
This map verifies the hypothesis of Corollary \ref{cB.1}, hence it is $0$. 
\end{proof}

As a variant, 

\begin{prop} We have an isomorphism 
\[P_{X\times Y}\iso R_{\pi_0(X)/F}(P_Y\times_F \pi_0(X))\times
R_{\pi_0(Y)/F}(P_X\times_F \pi_0(Y)).\]
\end{prop}

Since we are not going to use this, we leave the easy proof to the reader.

Proposition \ref{p4.2} a) endows $\SAbS$ with a symmetric monoidal 
structure compatible with its additive structure, hence also its full subcategory $\AbS$. From now on we concentrate on 
this latter category.

\begin{prop}\label{p4.3} The category $\AbS$ is symmetric 
monoidal (for $\otimes_\rep$) and pseudo-abelian. Its Kelly radical $\sR$
is monoidal and has square $0$. After tensoring with $\Q$, $\AbS/\sR$
becomes isomorphic to the semi-simple category product of the category of
abelian varieties up to isogenies and the category of $G_F$-$\Q$-lattices.
\end{prop}

Recall that the \emph{Kelly radical} $\sR$ of an additive catgegory $\sA$
is defined by 
\[\sR(A,B)=\{f\in \sA(A,B)\mid \forall g\in \sA(B,A)\; 1_A -gf\text{ is
invertible}\}\]
and that it is a [two-sided] ideal of $\sA$ \cite{kelly}.

\begin{proof} For the first claim, we just observe that kernels exist in the
category of commutative $F$-group schemes, and that a direct summand of an
abelian variety (\resp of a lattice) is an abelian variety (\resp a
lattice). For the second claim, consider the functor
\begin{align*}
T:\AbS&\to \Ab\times \Lat\\
\sA&\mapsto (\sA^0,\pi_0(\sA))
\end{align*}
where $\Ab$ and $\Lat$ are respectively the category of abelian
varieties and the category of lattices over $F$ (viewed, for example, as
full subcategories of the category of \'etale sheaves over $Sm/F$). This
functor is obviously essentially surjective. After tensoring with $\Q$,
it becomes full, because any extension 
\[0\to \sA^0\to \sA\to \pi_0(\sA)\to 0\]
is rationally split. Now the collection of sets 
\[\sI(\sA,\sB)=\{f:\sA\to\sB\mid T(f)=0\}\]
defines an ideal $\sI$ of $\AbS$. If
$f\in \sI(\sA,\sB)$, then $f$ induces a map 
\[\bar f:\pi_0(\sA)\to \sB^0\]
and this gives a description of $\sI$. From this description, it follows
immediately that $\sI^2=0$. In particular, $\sI\subseteq \sR$.

If we tensor with $\Q$, then $\Ab\times \Lat$ becomes semi-simple;
since $\AbS/\sI\otimes \Q$ is semi-simple and $\sI\otimes \Q$ is
nilpotent, it follows that $\sI\otimes\Q=\sR\otimes\Q$. In other words,
$\sR/\sI$ is torsion.

Let $f\in \sR(\sA,\sB)$. There exists $n>0$ such that $nf(\sA^0)=0$. But
$f(\sA^0)$ is an abelian subvariety of $\sB^0$, hence $f(\sA^0)=0$ and
$f\in \sI(\sA,\sB)$. So $\sR=\sI$.

If we endow the category $\Ab\times \Lat$ with the tensor structure
\[(A,L)\otimes (B,M)=(A\otimes M\oplus B\otimes L,L\otimes M)\]
then $T$ becomes a monoidal functor, which shows that $\sR=\sI$ is
monoidal. This completes the proof of Proposition \ref{p4.3}.\end{proof}

\begin{sloppypar}
\begin{rks} a) The morphisms in $\AbS$ are best represented in matrix
form:
\[\Hom(\sA,\sB)=
\begin{pmatrix}
\Hom(\sA_0,\sB_0)&\Hom(\pi_0(\sA),\sB_0)\\
0& \Hom(\pi_0(\sA),\pi_0(\sB))
\end{pmatrix}\]
(note that $\Hom(\sA_0,\pi_0(\sB))=0$). This clarifies the arguments in
the proof of Proposition \ref{p4.3} somewhat.

b) The Hom groups of $\Ab\times \Lat$ are finitely
generated $\Z$-modules. It follows from the proof of Proposition
\ref{p4.3} that, for $\sA,\sB\in \AbS$, $T(\Hom(\sA,\sB))$ has finite
index in $\Hom(T(\sA),T(\sB))$. In particular, for any $\sA\in\AbS$,
$\End(\sA)$ is an extension of an order in a semi-simple
$\Q$-algebra by an ideal of square $0$. 

c) The functor $T$ has the explicit section
\[(A,L)\mapsto A\oplus L.\]

This section is symmetric monoidal.
\end{rks}
\end{sloppypar}

\section{Chow birational motives and locally abelian schemes}\label{s6}

\subsection{The Albanese map}

For any smooth projective variety $X$, there is a canonical map
\begin{equation}\label{eq3}
\begin{CD}
CH_0(X)@>\Albv_X^F>> \sA_X(F).
\end{CD}
\end{equation}

Recall the construction of $\Albv_X$: the map $\phi_X$ of \eqref{eq6.2}
defines for any extension $E/F$ a map $X(E)\to \sA_X(E)$, still denoted by
$\phi_X$. When $E/F$ is finite, viewing $\sA_X$ as an \'etale
sheaf, we have a trace map
$Tr_{E/F}:\sA_X(E)\to \sA_X(F)$. Then $\Albv_X$ maps the class of a
closed point $x\in X$ with residue field $E$ to
$Tr_{E/F}\,\phi_{X}(x)$.

The map $\Albv_X$ is  injective for
$\dim X=1$ and surjective if
$F$ is algebraically closed. For a curve, this map corresponds to 
the isomorphism $\Pic_X\simeq \sA_X$, where $\Pic_X$ is the Picard scheme
of $X$; we then also have $\sA_X^0\simeq J_X$, where $J_X$ is the Jacobian
variety of $X$.

The functoriality of $\sA$ shows that there is a chain of isomorphisms
\begin{equation}\label{eq4}
\Phi_{X,Y}:\Hom(\sA_X,\sA_Y)\iso \Mor(X,\sA_Y)\iso \sA_Y(F(X))
\end{equation}
(the latter by Weil's theorem on extension of morphisms to abelian 
varieties \cite[Th. 3.1]{milne}), hence a canonical map
\begin{equation}\label{eq5}
\begin{CD}
CH_0(Y_{F(X)})@>\Albv_{X,Y}>> \Hom(\sA_X,\sA_Y)
\end{CD}
\end{equation}
which generalises \eqref{eq3}; more precisely, we have
\begin{equation}\label{eq7}
\Phi_{X,Y}\circ \Albv_{X,Y}=\Albv_Y^{F(X)}.
\end{equation}

On the other hand, there is an exact sequence
\begin{multline*}
0\to \sA_Y(\pi_0(X))=\Hom(\Z[\pi_0(X)],\sA_Y)\to \Hom(\sA_X,\sA_Y)\\
\to \Hom(\sA_X^0,\sA_Y)
\to \Ext^1(\Z[\pi_0(X)],\sA_Y)=H^1(\pi_0(X),\sA_Y)
\end{multline*}
and the map $\Hom(\sA_X^0,\sA_Y^0)\to \Hom(\sA_X^0,\sA_Y)$ is an
isomorphism. From this and \eqref{eq5} we get a zero sequence
\begin{equation}\label{eq6}
0\to CH_0(Y)\to CH_0(Y_{F(X)})\to \Hom(\sA_X^0,\sA_Y^0)\to 0.
\end{equation}

\begin{lemma}\label{l6.3} Let $Y,Z$ be two smooth projective varieties and
$\beta\in CH_0(Z_{F(Y)})$. Then the
following diagram commutes:
\[\begin{CD}
CH_0(Y)@>\beta_*>>CH_0(Z)\\
@V{\Albv_Y^F}VV @V{\Albv_Z^F}VV\\
\sA_Y(F)@>\Albv_{Y,Z}(\beta)_*>>\sA_Z(F).
\end{CD}\]
\end{lemma}

\begin{proof} Without loss of generality, we may assume that $\beta$ is given by
an integral subscheme $W$ in $Y
\times Z$. Then the composite $f=p_Yi_{W}$ is a proper surjective
generically finite morphism, where $p_Y$ denotes the projection and
$i_{W}$ is the inclusion of $W$ in $Y \times Z$. 

Let $V$ be an affine dense open subset of $Y$ such that $f_{|f^{-1}(V)}$
is finite. Any element of $CH_0(Y)$ may be represented by a zero-cycle
with support in $V$ (\cf \cite{roberts}), so it is enough to check
the commutativity of the diagram on zero-cycles on $Y$ of the form $y$,
where
$y\in V_{(0)}$. For such a $y$, we have $\beta_*y=p_*(f^{-1}(y))$, where
$p=p_Zi_W$. 

On the other hand, the composition $\Albv_{Y,Z}(\beta)_*\circ
(\Albv_Y^F)_{|V}$ may be described as follows: let $d$ be the degree of
$f_{|f^{-1}(V)}$, $f^{-1}(V)^{[d]}$ the $d$-fold symmetric power of
$f^{-1}(V)$ and $f^*:V\to f^{-1}(V)^{[d]}$ the map $x\mapsto f^{-1}(x)$.
Then
\[\Albv_{Y,Z}(\beta)_*\circ(\Albv_Y^F)_{|V}=\Sigma_d\circ
(\phi_Z)^{[d]}\circ p^{[d]}_*\circ f^*\]
where $\Sigma_d:\sA_Z^{[d]}\to \sA_Z$ is the summation map. The
commutativity of the diagram is now clear.\end{proof}

\subsection{The Albanese functor}

\begin{prop} \label{p6.1} The assignment $X\mapsto \sA_X$ defines via
\eqref{eq5} a symmetric monoidal
additive functor
\[\Albv:\Chow^\o\to \AbS\]
which becomes full and essentially surjective after tensoring with $\Q$.
\end{prop}

\begin{proof} Since $\AbS$ is pseudo-abelian, it suffices to construct the
functor on $\Cor^\o$. Let
$\alpha
\in CH_0(Y_{F(X)})$ and
$\beta
\in CH_0 (Z_{F(Y)})$. We want to show that
$\Albv_{X,Z}(\beta\circ\alpha)=\Albv_{Y,Z}(\beta)\circ
\Albv_{X,Y}(\alpha)$. But
$\beta$ induces a map 
\[\beta_*: CH_0(Y_{F(X)}) \to CH_0(Z_{F(X)}),\]
and we have the equality $\beta_*\alpha=\beta\circ\alpha$ (\cf proof
of Proposition \ref{l4.2'}). Hence, applying Lemma \ref{l6.3} in which we replace $F$
by $F(X)$, we get
\[\Albv^{F(X)}_{Z}(\beta\circ\alpha)=
\Albv^{F(X)}_{Z}(\beta_*\alpha)=\Albv_{Y,Z}(\beta)_*(\Albv^{F(X)}_Y(\alpha)).
\]

Applying now \eqref{eq7}, we get
\[
\Phi_{X,Z}\circ\Albv_{X,Z}(\beta\circ\alpha)
=\Albv_{Y,Z}(\beta)_*(\Phi_{X,Y}\circ\Albv_{X,Y}(\alpha)).
\]

On the other hand, the diagram
\[\begin{CD}
\sA_Y(F(X))@>\Albv_{Y,Z}(\beta)_*>> \sA_Z(F(X))\\
@A{\Phi_{X,Y}}A{\wr}A @A{\Phi_{X,Z}}A{\wr}A\\
\Hom(\sA_X,\sA_Y)@>\Albv_{Y,Z}(\beta)_*>> \Hom(\sA_X,\sA_Y)
\end{CD}\]
obviously commutes, which concludes the proof that $\Albv$ is a functor.
 
Compatibility with the monoidal structures follows from Proposition
\ref{p4.2} b). It remains to show the assertions on fullness and
essential surjectivity.

{\bf Fullness:} for any $Y$, the map $\Albv_Y^F\otimes \Q$ is surjective.
This follows from the case where $F$ is algebraically closed (in which
case
$\Albv_Y^F$ itself is surjective) by a transfer argument. Replacing the
ground field $F$ by $F(X)$ for some other $X$, we get that
$\Albv_{X,Y}\otimes \Q$ is surjective. This shows that the restriction
of $\Albv\otimes\Q$ to $\Cor^\o\otimes\Q$ is full; but the
pseudo-abelianisation of a full functor is evidently full (a direct
summand of a surjective homomorphism of abelian groups is surjective).

{\bf Essential surjectivity:} we first note that, after tensoring
with
$\Q$, the extension
\[0\to \sA^0\to \sA\to \pi_0(\sA)\to 0\]
becomes split for any $\sA\in \AbS$. Indeed the extension class
belongs to $\Ext^1_F(\pi_0(\sA),\sA^0)$; this group sits in an exact
sequence (coming from an Ext spectral sequence)
\begin{multline*}
0\to H^1(F,\Hom_{\bar F}(\pi_0(\sA)_{|\bar F},\sA^0_{|\bar
F}))\to \Ext^1_F(\pi_0(\sA),\sA^0)\\
\to H^0(F,\Ext^1_{\bar
F}(\pi_0(\sA)_{|\bar F},\sA^0_{|\bar F})).
\end{multline*}

Since the restriction $\pi_0(\sA)_{|\bar F}$ is a constant sheaf of free
finitely generated abelian groups, the group $\Ext^1_{\bar
F}(\pi_0(\sA)_{|\bar F},\sA^0_{|\bar F})$ is $0$, while the left group is
torsion as a Galois cohomology group. It is now sufficient to show
separately that $L$ and $A$ are in the essential image of $\Albv\otimes
\Q$, where $L$ (\resp $A$) is a lattice (\resp an abelian variety). 

A lattice $L$ corresponds to a continuous
integral representation $\rho$ of $G_F$. But it is well-known
that $\rho\otimes\Q$ is of the form $\theta\otimes \Q$, where $\theta$ is
a direct summand of a permutation representation of $G_F$. If $E$ is the
corresponding
\'etale algebra, we therefore have an isomorphism of $L$ with a direct
summand of $(\Albv\otimes\Q)(E)$.

Given an abelian variety $A$, we simply note that
\[A=\Albv(\tilde h(A))\]
where $\tilde h(A)$ is the reduced motive of $A$: $h(A)=\un\oplus \tilde
h(A)$, where the splitting is given by the rational point $0\in A(F)$.
\end{proof}

\begin{rk} Let $\sR$ be the Kelly radical of $\AbS$ (\cf
Proposition \ref{p4.3}). If
$F$ is a finitely generated field, the groups
$\sR(\sA,\sB)$ are finitely generated by the
Mordell-Weil-N\'eron theorem. To see this, note that if $L$ is a lattice
and $A$ an abelian variety, then
\[\Hom(L,A)\iso \Hom(L_{|\bar F},A_{|\bar F})^{G_F}\]
and that the right term may be rewritten as $B(F)$, where $B=L^*\otimes
A$ (compare Lemma \ref{l4.1}). Hence the Hom groups in $\AbS$ are
finitely generated as well. In this case, Proposition
\ref{p6.1} implies that, for any $M,N\in \Chow^\o$, \emph{the image of
the map
$\Albv_{M,N}$ has finite index in the group $\Hom(\Albv(M),\Albv(N))$}.
\end{rk}

\begin{lemma}\label{l9.1} Suppose that $Y$ is a curve. Then the map
\eqref{eq5} fits into an exact sequence
\begin{multline*}
0\to CH_0(Y_{F(X)})\by{\Albv_{X,Y}} \Hom(\sA_X,\sA_Y)\\
\to Br(F(X))\to Br(F(X\times Y))
\end{multline*}
where $Br$ denotes the Brauer group. In particular,
\eqref{eq5}$\otimes\Q$ is an isomorphism.
\end{lemma}

\begin{proof} First assume that $X$
is a point; then \eqref{eq5} reduces to \eqref{eq3}. Suppose first that
$F$ is separably closed. Then \eqref{eq3} is bijective (see comments at
the beginning of this section). In the general case, let $F_s$ be a
separable closure of $F$, and $G=Gal(F_s/F)$. Since
$\sA_Y$ is a sheaf for the \'etale topology, we get a commutative diagram
\[\begin{CD}
CH_0(Y_s)^G@>\Albv_Y^{F_s}>\sim> \sA_Y(F_s)^G\\
@A{}AA @A{\wr}AA\\
CH_0(Y)@>\Albv_Y^F>> \sA_Y(F)
\end{CD}\]
where $Y_s=Y\times_F F_s$ and the top horizontal and right vertical maps
are bijective. The lemma then follows from the classical exact sequence
\[0\to CH_0(Y)\to CH_0(Y_s)^G\to Br(F)\to
Br(F(Y)).\]

The case where $X$ is not necessarily a point now follows from this special case and the construction of \eqref{eq5}.
\end{proof}

\begin{thm}\label{p6.2} Let $\Chow^\o_{\le 1}$ denote the thick subcategory
of 
\break $\Chow^\o$ generated by motives of varieties of dimension $\le 1$, and
let $\iota:\Chow^\o_{\le 1}\allowbreak \to \Chow^\o$ be the canonical inclusion.
Then\\ 
a) After tensoring morphisms with $\Q$, $\Albv\circ \iota:\Chow^\o_{\le
1}\to \AbS$ becomes an equivalence of categories.\\ 
b) Let $j$ be a quasi-inverse. Then $\iota\circ j$ is right adjoint to
$\Albv$.
\end{thm}

\begin{proof} a) The full faithfulness follows from Lemma \ref{l9.1}. For the
essential surjectivity, we may reduce as in the proof of Proposition
\ref{p6.1} to proving that lattices and abelian varieties are in the
essential image. For lattices, this is proven in \loccit. For an abelian
variety $A$, use the fact that $A$ is isogenous to a quotient
of the Jacobian of a curve, and Poincar\'e's complete reducibility
theorem.

b) Let $(M,\sA)\in \Chow^\o_{\le 1}(F,\Q)\times \AbS(F,\Q)$. To produce a
natural isomorphism $\Chow^\o_{\le 1}(F,\Q)(M,\iota j(\sA))\simeq
\AbS(F)(\Albv(M),\sA)\otimes \Q$, it is sufficient by a) to handle the case
$M= h^\o(X), \sA=\sA_Y$ for some smooth projective curves $X,Y$.
Then the isomorphism follows from the isomorphisms 
\eqref{eq4} and from Lemma
\ref{l9.1}.\end{proof}

\begin{rks} a) Of course the functor $\iota\circ j$ is not a tensor
functor (since its image is not closed under tensor product).\\
b) In particular, the inclusion functor $\iota$ has the left adjoint
 $j\circ \Albv$. This is a birational version of Murre's
results for effective Chow motives (\cite{murre}, \cite[\S 2.1]{murre2},
see also \cite[\S 4]{scholl}). Beware however that we have taken the opposite to
usual convention for the variance of Chow motives (our functor
$X\mapsto h(X)$ is covariant rather than contravariant), so the direction
of arrows has to be reversed with respect to Murre's work.
\end{rks}

\appendix 
\section{Complements on localisation of categories}

\subsection{Localisation of symmetric monoidal categories}

\begin{lemma}\label{lA.2} a) Localisation commutes with products of
  categories for sets of morphisms containing all identities\footnote{We thank M. Bondarko for pointing out the importance of the identities.}.\\ 
b) Let $T_0, T_1:\sC\rr \sD$ be two functors and $f:T_0\Rightarrow
T_1$ a natural  
transformation. Let $S, S'$ be collections of morphisms in $\sC$ and
$\sD$ such that $T_i(S) \subseteq S'$, so that $T_0$ and $T_1$ pass to
localisation. Then $f$ remains a  
natural transformation between the localised functors.
\end{lemma}

\begin{proof} a)  Let $S_i$ be a collection of morphisms in $\sC_i$ for $i=1,2$, such that $S_i$ contains the identities of all objects of $\sC_i$. Then $S_1\times S_2$ is generated by $S_1$ and $S_2$ in the sense that the equality
\[(s_1,s_2)=(s_1,1) \circ (1,s_2)\]
holds in  $S_1\times S_2$ for any pair  $(s_1,s_2)$. The conclusion easily follows
 (\cf \cite[Lemma 2.1.7]{maltsin}). b) is true because $f$ commuted with the members
of $S$, hence it now  
commutes with their inverses.\end{proof}

\enlargethispage*{30pt}

\begin{prop}\label{pA.2} Let $\sC$ be a category with a product
  $\bullet: \sC\times \sC\to \sC$, and let $S$ be a  
collection of morphisms in $\sC$ containing all identities. Assume that $S\bullet S\subseteq S$. Then\\
a) There is a unique product  $S^{-1}\sC\times S^{-1}\sC\to S^{-1}\sC$
such that the localisation functor $P_S:\sC\to S^{-1}\sC$ commutes
with the two products.\\ 
b) If $\bullet$ is monoidal (\resp braided, symmetric,  unital),  the
induced product on $S^{-1}\sC$ enjoys the same properties and $P_S$ is
monoidal (\resp braided, symmetric, unital). 
\end{prop}

\begin{proof} a) follows from Lemma \ref{lA.2} a); b)  follows from Lemma
\ref{lA.2} b).\end{proof}

\subsection{Semi-additive categories}

This subsection is a reformulation of \cite[Ch. VIII, \S 2]{maclane}, see also \cite[\S 18 and beginning
  of \S 19]{maclane1}.

\begin{lemma}\label{lA.1} a) For a category $\sA$, the following conditions are equivalent:
\begin{thlist}
\item $\sA$ has a $0$ object (initial and final), binary products and coproducts, and for any $A,B\in\sA$, the map
\[A\coprod B\to A\times B\]
given on $A$ by $(1_A,0)$ and on $B$ by $(0,1_B)$ is an isomorphism.
\item $\sA$ has finite products, and for any $A,B\in \sA$, $\sA(A,B)$ has a structure of a commutative monoid, and 
composition is distributive with respect to these monoid laws.
\item Same as {\rm (ii)}, replacing product by coproduct.  
\end{thlist}
We then say that $\sA$ is a \emph{semi-additive category} and write $A\oplus B$ for the product or coproduct of two 
objects $A,B$.\\
b) If $\sA$ is a semi-additive category, the law $(A,B)\mapsto A\oplus B$ endows $\sA$ with a canonical unital symmetric 
monoidal structure.
\end{lemma}

\begin{proof} a) By duality, we only need to show (i) $\iff$ (ii). (i) $\Rightarrow$ (ii) follows from \cite[Ch. VIII, \S 
2, ex. 4 (a)]{maclane}: recall that for two morphisms $f,g:A\to B$ in $\sA$, Mac Lane defines their sum $f+g$ as the 
composition
\[
\xymatrix{
A\ar[r]^{\Delta_A}\ar[dd]_{f+g}&A\times A\ar[dr]_{f\times g}\\
&&B\times B\\
B&B\coprod B\ar[l]_{\nabla_B}\ar[ur]^\sim
}
\]
where $\Delta_A$ is the diagonal and $\nabla_B$ is the codiagonal.

As for (ii) $\Rightarrow$ (i), it is implicit in the proof of \cite[Ch. VIII, \S 2, Th. 
2]{maclane}. Indeed, Mac Lane defines a biproduct of two objects $A,B\in \sA$ as a diagram
\[A
\begin{smallmatrix}
p_1\\
\leftrightarrows\\
i_1
\end{smallmatrix} 
C 
\begin{smallmatrix}
p_2\\
\leftrightarrows\\
i_2
\end{smallmatrix} 
B\]
satisfiyng $p_1i_1=1_A$,  $p_2i_2=1_B$ and $i_1p_1+i_2p_2=1_C$. Let us say that such a diagram is a \emph{biproduct*} if 
the further identities $p_1i_2=0$ and $p_2i_1=0$ hold. Then, Mac Lane proves that a biproduct* is a product and that a 
product 
is a biproduct*. Dually, a biproduct* is the same as a coproduct, hence binary products and coproducts are canonically 
isomorphic, and one checks from his proof that the isomorphism is given by the map of (i).

(Let us clarify that Mac Lane proves that a biproduct is a biproduct* if the addition law on morphisms has the 
cancellation property; but we don't use this part of his proof.)

b) This is obvious: already finite products or coproducts define a canonical symmetric monoidal structure.
\end{proof}

Define a \emph{semi-additive functor} between two semi-additive categories $\sA,\sB$ as a functor $F:\sA\to \sB$ which 
preserves addition of morphisms. Note that any semi-additive functor preserves $\oplus$, by the characterisation of 
biproducts via equations (see proof of Lemma \ref{lA.1} a)). 

\subsection{Localisation of $R$-linear categories}

\begin{thm}\label{tA.1} Let $\sA$ be a semi-additive category and $S$ a family
of morphisms of $\sA$, containing all identities and stable under $\oplus$. Then  $S^{-1}\sA$ and the localisation functor $P_S:\sA\to S^{-1}\sA$ are semi-additive.
\end{thm}

\begin{proof} We use the characterisation (i) of semi-additive categories in Lemma \ref{lA.1}: by \cite[1.3.6 
and 2.1.8]{maltsin}, $P_S$ preserves products and coproducts, and transforms the isomorphisms $A\coprod B\iso A\times B$ 
into isomorphisms.
\end{proof}

To ``catch" additive categories (as opposed to semi-additive categories), we
could do as in Mac Lane \cite{maclane1} and postulate the existence of an
endomorphism $-1_A$ for each object $A$. We prefer to do this more generally by
dealing with $R$-linear categories, where $R$ is an arbitrary ring (an
$R$-linear category is simply a semi-additive $R$-category).

More precisely, let $\sA$ be an $R$-linear category. Then in particular:

\begin{itemize}
\item $\sA$ is a semi-additive category.
\item It enjoys an action of the multiplicative monoid underlying $R$, \ie there is a homomorphism of monoids $R\to 
\End(Id_A)$, where $\End(Id_A)$ is the monoid of natural transformations of the identity functor of $A$.
\item For $\lambda\in R$ and $A\in \sA$, let $\lambda_A$ denote the
corresponding endomorphism of $A$.
Then we have identities
\begin{equation}\label{add}
(\lambda+\mu)_A=\lambda_A+\mu_A.
\end{equation}
\end{itemize}

Conversely, the following lemma is straightforward.

\begin{lemma} Let $\sA$ be a semi-additive category provided with an action of $R$
verifying \eqref{add}. Then $\sA$ is an $R$-linear category.\qed
\end{lemma}

From this lemma, it follows:

\begin{thm} \label{tA.2} Theorem \ref{tA.1} extends to
$R$-linear categories.\qed
\end{thm}

\enlargethispage*{20pt}

\subsection{Localisation and pseudo-abelian envelope}

\begin{lemma}\label{lA6.1} Let $\sA$ an additive category and $S$ a family of morphisms in
$\sA$, stable under direct sums. Let $\sA\to \sA^\natural$ denote the pseudo-abelian envelope
of $\sA$, and let us denote by $S^\natural$ the set of direct summands of members of $S$ in
$\sA^\natural$. Then the natural functors
\[(S^{-1}\sA)^\natural\to(S^{-1}(\sA^\natural))^\natural\to((S^\natural)^{-1}(\sA^\natural))^\natural\]
are equivalence of categories.
\end{lemma}

\begin{proof} All categories are universal for additive functors $T$ from $\sA$
to a pseudo-abelian category such that $T(S)$ is invertible.\end{proof}

\subsection{Localisation and group completion}

\begin{lemma} Let $\sA$ be a semi-additive category. There exists an additive category $\sA^+$ and a semi-additive functor $\iota:\sA\to \sA^+$ with the following $2$-universal property: any semi-additive functor from $\sA$ to an additive category factors through $\iota$ up to a unique natural isomorphism.\\
A model of $\sA^+$ may be given as follows: the objects of $\sA^+$ are those of $\sA$; if $A,B\in \sA$, then $\sA^+(A,B)$ is the group completion of the commutative monoid $\sA(A,B)$.\\
The category $\sA^+$ is called the \emph{group completion} of $\sA$.
\end{lemma}

The proof is straighforward and omitted.

\begin{prop} Let $\sA$ be a semi-additive category, and let $S$ be a family of morphisms in $\sA$, containing the identities and stable under direct sums. Keep writing $S$ for the image of $S$ in the group completion $\sA^+$. Then the functor $S^{-1}\iota:S^{-1} \sA \to S^{-1} (\sA^+)$ induces an equivalence of categories
\[\tilde \iota:(S^{-1} \sA)^+\iso S^{-1} (\sA^+).\]
Here we use the structure of semi-additive category on $S^{-1}\sA$ given in Theorem \ref{tA.1}.
\end{prop}

\begin{proof} The existence of $\tilde\iota$ follows from the universal property of group completion. A quasi-inverse to $\tilde\iota$ is obtained by group-completing the functor $\sA\to S^{-1}\sA$ (which is semi-additive by Theorem \ref{tA.1}), and then extending the resulting functor to $S^{-1} (\sA^+)$.
\end{proof}

\end{document}